\documentclass[11pt]{epiarticle}
\usepackage{epimath}

\setpapertype{A4}

\usepackage{hyperref}
\hypersetup{linkcolor=blue,citecolor=red,filecolor=green,urlcolor=magenta} 

\usepackage[english,french]{babel}
\addto\captionsfrench{%
   \def\figurename{Figure}%
   \def\tablename{Tableau}%
}
\usepackage[T1]{fontenc}
\usepackage[utf8]{inputenc}
\usepackage{caption}

\usepackage[pdftex]{graphicx}  

\usepackage{apacite}

\usepackage{tabularx}
\newcolumntype{C}{>{\centering}X}

\usepackage[titletoc]{appendix}

\newcommand{\R}{\mathbb{R}}
\newcommand{\N}{\mathbb{N}}
\newcommand{\m}{\textperiodcentered}
\newtheorem{exo}{Exercice}
\usepackage{comment}
\usepackage{csquotes}
\usepackage{natbib}
\usepackage{enumitem}
\usepackage{array}
\usepackage{arydshln}
\usepackage{csquotes}

\title{Évaluation entre pairs en mathématiques : activité d'étudiant\m es lors d'évaluations de preuves}
\titlemark{Évaluation entre pairs en mathématiques}
\author{Juliette Veuillez-Mainard \&\ Simon Modeste}
\authormark{J.\ Veuillez-Mainard \&\ S.\ Modeste}
\date{} 
\journal{\'EpiDEMES  -- (aaaa), ---} 
\acceptation{soumis le 21/10/2023, accepté le jj/mm/aaaa, version définitive le jj/mm/aaaa}

\begin{document}

\selectlanguage{french}

\maketitle



\thanks{Ce travail a été financé en partie par l'Agence Nationale de la Recherche (ANR) grâce au projet CaMeLOt ANR-20-CHIA-0001-01.}

\begin{prelims}

\selectlanguage{english}

\def\abstractname{Abstract}
\abstract{This article deals with peer-assessment in the context of higher education teaching in mathematics, and examines the nature of student activity when assessing work produced by peers. After an overview of research on peer assessment, we propose an experiment with students in a specific post-secondary scientific class, and analyze the activity of the students involved. Our analyses are based on a priori analyses of the proposed tasks and tools from research on students personal mathematical work. We base on written notes of assessments, observations of pairs in assessment situations, and answers to a questionnaire on the perception of their activity.  We have identified some specific student activities during the peer assessment process, in particular related to the analysis and correction of proofs, but also to the awareness of the issues associated with evaluation. This enables us to argue for the potential of peer assessment in higher mathematics education.}

\keywords{Didactics of mathematics, Peer assessment, Proof and proving, Students' activity, University education, Real numbers}

\medskip

\selectlanguage{french}
\def\abstractname{Résumé}
\abstract{Cet article traite de l’évaluation entre pairs (\textit{peer-assesment}) dans le contexte de l’enseignement supérieur en mathématiques, et étudie la nature de l’activité d’étudiant\m es lors de l’évaluation de copies produites par des pairs. Après un tour d’horizon des recherches sur l’évaluation par les pairs, nous proposons une expérimentation en classe préparatoire MPSI et nous analysons l’activité des étudiant\m es impliqué\m es. Nos analyses s’appuient sur des analyses a priori des tâches proposées et des outils issus de recherches sur le travail personnel d’élèves en mathématiques. Nous nous basons sur les traces écrites des évaluations, des observations de binômes en situation d’évaluation, et les réponses à un questionnaire sur la perception de leur activité.  Nous identifions certaines activités spécifiques des étudiant\m es lors du processus d’évaluation de pairs, notamment liées à l’analyse et la correction de preuve, mais aussi liées à la prise de conscience des enjeux associés à l'évaluation. Ceci nous permet d’argumenter du potentiel de l’évaluation entre pairs pour l’enseignement des mathématiques dans le supérieur.}

\def\keywordsname{Mots-Clés}
\keywords{Didactique des mathématiques, Évaluation par les pairs, Preuve et démonstration, Activité d'étudiants, Enseignement supérieur, Nombres réels}

\selectlanguage{french}

\tableofcontents

\end{prelims}


\section{Introduction}

Les apprentissages des élèves et étudiant\m es\footnote{Nous avons choisi d'utiliser le terme générique d'\textit{élève}, à entendre comme incluant l'enseignement supérieur, dans toute la partie de généralités sur l'évaluation entre pairs, et le terme d'\textit{étudiant\m es} pour nos expérimentations qui portent sur le début du supérieur.} sont évalués dès l'enseignement primaire et tout au long de la scolarité. Ces évaluations peuvent prendre différentes formes en fonction des types de connaissances que l'on souhaite évaluer, de l'âge des élèves, ou encore de leur parcours. Elles sont un outil de mesure des connaissances et des apprentissages des élèves, et peuvent servir au développement et à l'adaptation de l'enseignement. Plusieurs types d'évaluations peuvent être mises en place en fonction des objectifs de l'évaluation, par exemple des devoirs en classe ou à la maison, ou des présentations orales ou affichées. La question de l'implication des élèves dans l'évaluation a été discutée dans la recherche, et diverses études ont montré l'intérêt des processus d'auto-évaluation et d'évaluation par les pairs, notamment leur fiabilité et la diversité des apprentissages qu'ils peuvent favoriser.

Nous nous intéresserons dans cet article à l'évaluation par les pairs dans l'enseignement supérieur. D'après la définition de \citet{topping_self_2003}, \enquote{l'évaluation par les pairs est un dispositif dans lequel les apprenants et/ou travailleurs sont amenés à examiner et spécifier la valeur ou la qualité d'un produit ou d'une performance d'autres apprenant\m s et/ou travailleur\m ses de même statut} \citep[p.~65, notre traduction]{topping_self_2003}. Dans le cadre d'évaluation par les pairs dans un contexte d'enseignement, les élèves sont donc amenés à évaluer d'autres élèves de même niveau. Les précédentes recherches sur ce sujet nous invitent à supposer que cette évaluation est un soutien à l'apprentissage. Nous faisons l'hypothèse qu'être mis en situation d'évaluation de pairs change l'activité des élèves, et nous nous demandons donc quelle activité les étudiant\m es mettent en place pendant ce processus, et quels apprentissages cela pourrait alors favoriser.

Nous présentons dans la première partie de cet article un tour d'horizon succinct de la recherche sur l'évaluation par les pairs, complété par une classification des \enquote{processus de vérification en mathématiques chez les élèves} issue de \citet{coppe_processus_1993}. Ceci permettra de formuler nos questions de recherche. Dans la deuxième partie, nous présentons l'expérimentation que nous avons menée et précisons sa méthodologie. Les analyses \textit{a priori} du problème et des évaluations de ce problème sont ensuite développées. Nous exposons dans la troisième partie les analyses de l'expérimentation, avant de présenter les résultats du questionnaire dans la quatrième partie. En conclusion, nous revenons sur les principaux résultats de notre étude, les perspectives de recherche, et des idées et conseils pour la mise en place d'évaluations par les pairs.

\section{Recherches sur l'évaluation par les pairs et positionnement}

La notion d'évaluation peut être définie de manière générale comme le propose \citet{de_ketele_en_1986} :

\begin{quote}
    L'évaluation est le processus qui consiste à recueillir un ensemble d'informations pertinentes, valides et fiables, puis à examiner le degré d'adéquation entre cet ensemble d'informations et un ensemble de critères choisis adéquatement en vue de fonder la prise de décision. \citep[p.~266]{de_ketele_en_1986}
\end{quote}


Il est courant de distinguer l'\emph{évaluation sommative}, qui a une fonction de certification et de contrôle, de l'\emph{évaluation formative}, qui vise à être un soutien explicite aux apprentissages \citep{morrissette_panorama_2014}. Le concept d'évaluation formative a été introduit en 1967 par Scriven, et décrit une évaluation mise en place dans le but de fournir des données qui vont permettre des adaptations de l'enseignement \citep{allal_formative_2005}. Enfin, lorsque le but premier de l'évaluation est l'apprentissage des élèves plutôt que la mesure de leurs connaissances, on parle d'\textit{évaluation soutien d'apprentissage}. Plus précisément, lors de la troisième conférence internationale sur l'Évaluation soutien d'Apprentissage, cette évaluation a été définie de la manière suivante :

\begin{quote}
    L'\'Evaluation soutien d'Apprentissage fait partie des pratiques quotidiennes des élèves et enseignants qui, individuellement ou en interaction, recherchent, réfléchissent sur et réagissent à l'information provenant d'échanges, démonstrations et observations afin de favoriser les apprentissages en cours. \citep[p.~102]{allal_assessment_2009}
\end{quote}

De nombreuses recherches sur l'évaluation formative ont été réalisées, ainsi que plusieurs revues et panorama de ces recherches. On pourra par exemple consulter les synthèses de \citet{black_assessment_1998}, \citet{allal_formative_2005} et \citet{morrissette_panorama_2014}, sur lesquelles nous nous appuyons dans la suite, et qui intègrent des travaux allant du primaire au supérieur, toutes disciplines confondues, incluant les mathématiques.

\subsection{Généralités sur l'évaluation par les pairs}
\label{generalites}

L'évaluation par les pairs est une forme d'évaluation dans laquelle les élèves sont amenés à évaluer les productions d'autres élèves de même niveau. C'est un processus distinct de l'enseignement par les pairs, dans lequel les élèves ne sont pas mis en position d'évaluateur.

Plusieurs éléments peuvent varier dans la mise en place d'une évaluation par les pairs : nature des apprentissages évalués, nature des supports évalués... Nous nous concentrerons ici sur l'évaluation de travaux individuels dans le cadre de l'enseignement. Dans la suite de cette section, nous nous appuyons sur les méta-analyses et revues de littérature de \citet{falchikov_student_2000} et \citet{topping_peer_1998, topping_self_2003}. En cohérence avec la littérature sur le sujet, nous distinguerons l'évaluation par les pairs (peer-assesment) et la notation par les pairs (peer-grading). Deux enjeux principaux se présentent : la \textit{validité} et la \textit{fiabilité} de l'évaluation, et ses apports possibles à l'apprentissage.

\paragraph{Fiabilité et validité}

Une partie de la littérature s'interroge sur la capacité des élèves à évaluer et noter leurs pairs. La notion de \emph{validité} réfère à la cohérence de l'évaluation (souvent la notation) des élèves par rapport à celle d'une référence experte (un\m e enseignant\m e par exemple). La notion de \emph{fiabilité}, en fait assez peu étudiée dans la littérature, réfère à la variabilité des évaluations (souvent des notations) entre élèves-évaluateurs, dans des conditions similaires.

Dans le cas de notations, l'étude de la validité s'appuie souvent sur des indicateurs statistiques, notamment la corrélation ou l'effect size. Ce dernier caractérise l'ampleur de la différence de note moyenne entre les élèves et les enseignants sur une tâche particulière \citep[p.~6]{sadler_impact_2006}. Les recherches concluent que, sous certaines conditions, les élèves peuvent produire des évaluations valides. Elles identifient des critères influant significativement sur la validité, par exemple la nature du produit évalué ou le degré d'avancement du cours (la significativité de ce second critère est discutée par \citet{falchikov_student_2000}). Il semble aussi que la familiarité des élèves avec les critères de notation tend à améliorer la validité de l'évaluation.

\paragraph{Bénéfices pour l'apprentissage}

Le processus d'évaluation par les pairs semble favoriser l'apprentissage des élèves. On observe d'une part des progrès du niveau effectif des élèves sur une tâche après la mise en place d'une évaluation par les pairs. Il semble aussi que ce dispositif peut impliquer de nouveaux questionnements chez les élèves, et permet de détecter plus tôt des conceptions erronées. De plus, elle implique en général que les élèves passent plus de temps à pratiquer une tâche.

D'autre part, l'évaluation par les pairs permet des gains métacognitifs pour les élèves. La présence d'un retour (feedback) personnalisé permet d'améliorer la capacité d'autoévaluation. En complément d'une évaluation classique, l'évaluation par les pairs peut apporter des feedbacks individualisés plus nombreux et réguliers, ce qui peut être un gain pour l'élève évalué\m e. Ce processus permet enfin d'observer des gains sociaux et de communication.

\subsection[Un regard en didactique des mathématiques sur l'évaluation par les pairs]{Un regard en didactique des mathématiques sur l'évaluation par les pairs : travaux de  Coppé et Moulin}

Les résultats mentionnés jusqu'à présent, même s'ils intègrent des études en mathématiques, ne sont pas spécifiques au contexte de la classe de mathématiques, et encore moins à la nature des connaissances mathématiques qui pourraient être évaluées.

En fait, peu de travaux en didactique des mathématiques se sont intéressés à l'évaluation entre pairs. Parmi eux, intéressons nous aux recherches réalisées par \citet{coppe_evaluation_2017}. Les autrices ont mis en place une séquence d'évaluation par les pairs et de débats argumentatifs dans une classe de sixième. Après avoir cherché un problème, les élèves se positionnent sur la validité des réponses des autres élèves. La classe est ensuite amenée à discuter collectivement des réponses affichées.

Le problème sélectionné porte sur les aires, et fait intervenir des connaissances sur les fractions dans sa résolution et dans les vérifications. Les autrices ont montré que le dispositif permet a priori de \enquote{provoquer [\dots] la mobilisation de connaissances autres que celles utilisées pour résoudre le problème} \citep[p.~326]{coppe_evaluation_2017}. Cependant, le problème s'est avéré trop complexe pour certains élèves, qui n'ont pas réussi à évaluer des démarches différentes de la leur. 


Dans l'analyse des processus de vérification mis en place par les élèves au moment de la correction, les autrices s'appuient sur le travail de thèse de \citet{coppe_processus_1993} définissant une \emph{typologique des processus de vérification}. Deux types sont distingués :

\begin{quote}
    Nous appellerons processus de vérification interne tout processus de vérification mettant en jeu des savoirs ou des savoir-faire typiquement mathématiques, ne dépendant pas nécessairement de la situation dans laquelle on les utilise. [\dots] Les autres processus, qui utilisent des connaissances portant sur d'autres savoirs ou savoir-faire moins mathématiques (notamment ceux qui n'utilisent pas seulement la logique du problème mais qui dépendent davantage du contrat) seront appelés des vérifications de type externe. Ces processus ont la propriété d'être généralisables à tous les problèmes ou à des classes de problèmes très étendues [\dots]. Certains peuvent être très courts et/ou se limiter à des arguments simples. \citep[p.~36]{coppe_processus_1993}
\end{quote}

Des exemples de processus de vérification externes sont les vérifications sociales dues au travail de groupe, ou les vérifications faisant référence au texte de l'exercice, comme le type d'unité pour une aire. En effet, ces processus \enquote{ne nécessitent pas \textit{a priori} d'entrer dans la stratégie de résolution du problème.} \citep[p.~12]{coppe_evaluation_2017}. D'un autre coté, on retrouve les arguments de vérification internes tels que les vérifications techniques du résultat ou les vérifications par changement de cadre.

La typologie proposée fournira un premier axe d'analyse pour l'expérimentation que nous présentons dans la suite de ce texte.

\subsection{Bilan et formulation des questions de recherche}

La littérature indique que les élèves, sous certaines conditions, fournissent des évaluations de pairs valides. Elle met à jour certains critères qui influent significativement sur la validité des évaluations et sur les apprentissages qui peuvent se produire. L'évaluation par les pairs permet par exemple aux élèves d'améliorer leur performance sur une tâche particulière, ou de développer des capacités métacognitives comme l'autoévaluation. Ce type d'évaluation semble donc orienté vers des fins formatives, pour les élèves évalués (rétroactions sur leur travail) mais aussi pour les élèves évaluateurs, ce qui fait la spécificité de l'évaluation par les pairs. Notons aussi que plusieurs articles donnent des conseils pratiques quant à la mise en place d'une telle évaluation (nous y reviendrons en fin d'article).



De nombreuses recherches se concentrent donc sur la validité et la fiabilité de l'évaluation entre pairs, à différents niveaux et dans différentes disciplines. Nous nous focalisons sur l'évaluation entre pairs en mathématiques, plus spécifiquement sur l'activité des étudiant\m es au cours de ce processus, un sujet qui semble peu abordé dans la littérature.

Nous faisons l'hypothèse que l'activité des étudiant\m es est différente lors d'une évaluation de pairs par rapport à la réalisation de la tâche évaluée elle-même, en particulier dans le travail de preuve. Nous nous concentrons dans cet article sur l'évaluation par les pairs, en mathématiques, dans les première années de l'enseignement supérieur. Ceci nous amène à poser les questions de recherche suivantes :

\begin{itemize}
    \item Quelle est l'activité d'étudiant\m es en situation d'évaluer des pairs ? Et en quoi diffère-t-elle de celle dans des situations plus courantes ?
    \item Quels éléments les étudiant\m es prennent-il\m elles en compte dans leurs évaluations ?
    \item Quelles corrections les étudiant\m es apportent-il\m elles à une copie qu'il\m elles évaluent ?
    \item À quelles ressources font appel les étudiant\m es en situation d'évaluer des pairs ?
    \item De quelle manière la position de correcteur\m rice contribue-t-elle à l'apprentissage des étudiant\m es ?
\end{itemize}

Pour apporter des éléments de réponse à ces questions, nous avons expérimenté un processus d'évaluation entre pairs en première année d'enseignement supérieur scientifique (plus précisément en CPGE MPSI\footnote{Classe Préparatoire aux Grandes Écoles, filière Mathématiques-Physique-Sciences de l'Ingénieur}), et avons recueillis différentes données : des résolutions d'exercices et les évaluations de ces exercices par les étudiant\m es, des vidéos d'entretiens d'explicitation de binômes d'étudiant\m es en situation d'évaluation de pairs, ainsi que les réponses des étudiant\m es à un questionnaire interrogeant la perception des étudiant\m es de l'activité d'évaluation de pairs.


\section{Expérimentation et méthodologie}

Cette section est consacrée à la présentation de l'expérimentation que nous avons menée. Nous commençons par expliciter notre positionnement théorique.

\subsection{Cadrage théorique}

Notre projet est d'étudier l'activité d'étudiant\m es mis\m es en situation d'évaluation de pairs. Notre hypothèse, en cohérence avec les recherches existantes, est que mettre les étudiant\m es dans une telle situation a des effets sur la nature de leurs activités, et notamment engendre des activités non-usuelles (sur les plans mathématique et cognitif, mais aussi méta-mathématique et méta-cognitif) que nous souhaitons identifier.
Pour cela, nous nous positionnons implicitement, sans nécessairement en reprendre les outils, dans le paradigme de la théorie de l'activité, telle qu'elle peut-être utilisée en didactique des mathématiques \cite{robert_activites_2017}. Ainsi, si nous adoptons principalement un modèle d'élève sujet épistémique, nous ne nous interdisons pas d'intégrer d'autres dimensions du sujet (psychologiques, sociales...) dans notre analyse, bien que nous ne cherchons pas à les catégoriser.

Nous nous appuierons sur une méthodologie de type ingénierie didactique \cite{artigue_ingenierie_1988}, qui permet un contrôle sur la conception et l'expérimentation de la tâche proposée, offre dans un mode de validation interne basé sur la comparaison analyse a priori - analyse a posteriori. Ceci nous permettra aussi d'étudier les activités possibles et activités réelles des étudiant\m es face à la tâche proposée. Nous intégrerons à ces analyses la typologie de \cite{coppe_processus_1993}, déjà utilisée par \cite{coppe_evaluation_2017} pour étudier des processus d'évaluation entre pairs, que nous adapterons à notre situation.

\subsection{Dispositif expérimental}

Le public sélectionné est composé d'étudiant\m es en début de l'enseignement supérieur, qui sont en apprentissage d'une pratique avancée des mathématiques. L'expérimentation s'est déroulée dans une classe préparatoire MPSI constituée de 42 étudiant\m es. Nous sommes conscients qu'un contrat didactique d'évaluation bien spécifique existe en CPGE (évaluations écrites régulières, attention forte à la rédaction, préparation aux "normes d'évaluation" des concours, rôle formatif de l'évaluation\dots), et nous en discuterons en conclusion. Cependant ce contexte peut aussi être pertinent pour notre étude, puisqu'il s'agit d'identifier l'activité que peuvent développer des étudiant\m es dans un contexte où il\m elles ont construit un rapport à l'évaluation propice à entrer dans un processus d'évaluation entre pairs.

Pour notre expérimentation, nous avons proposés des exercices aux étudiant\m es, puis les avons mis\m es en situation d'évaluation entre pairs. L'expérimentation en elle-même s'est déroulée pendant trois séances d'une heure de cours, complétées de moments d'évaluation à la maison pour la majorité des étudiant\m es. Les séances se sont déroulées de la manière suivante.
 
\begin{itemize}
    \item Pendant la première séance, l'expérimentation a été présentée aux étudiant\m es et il\m elles ont cherché et rédigé les exercices proposés. Nous avons ensuite récupéré et numérisé les copies.
    \item Au cours de la deuxième séance une discussion collective a été menée, d'une part sur la correction des exercices, et d'autre part sur les principes généraux de l'évaluation.
    \item Durant la semaine suivante, les étudiant\m es ont tous\m tes reçu deux copies à évaluer en travail à la maison, et certain\m es étudiant\m es volontaires ont réalisé leurs évaluations en binôme sous observation. Toutes les copies corrigées ont été récupérées et nous avons complété l'évaluation pour que les corrections rendues soient valides et complètes.
    \item La troisième et dernière séance a été consacrée à une correction complète des exercices, et à des commentaires sur les erreurs récurrentes rencontrées dans les copies.
\end{itemize}

Pour guider les étudiant\m es dans l'évaluation, en amont de la séance 2, nous avons réalisé une feuille d'"éléments de corrections", disponible en annexe. Elle a été conçue comme un support de réflexion pour l'évaluation et présente différentes démarches pour résoudre les exercices, ainsi que quelques éléments théoriques. Pendant l'évaluation, tous\m tes les étudiant\m es ont également eu accès à leurs notes de cours. À l'issu de la séance 3, nous avons transmis aux étudiant\m es un document présentant différentes résolutions possibles des exercices proposés, ainsi qu'une compilation des erreurs les plus fréquemment rencontrées.

Nous détaillons dans la suite les données recueillies, les méthodes d'analyse, les exercices proposés, ainsi que les analyses a priori.

\paragraph{Données recueillies}

À l'issue de la première séance, nous avons numérisé toutes les copies produites par les étudiant\m es. Sur le plan du dispositif, cela permet de faire évaluer une même copie par plusieurs étudiant\m es, et de faire évaluer plusieurs copies à chacun\m e. Sur le plan de la recherche, cela nous permet d'accéder au travail initial produit par chaque étudiant\m e.

À l'issue de la deuxième séance, nous avons donné aléatoirement à chaque étudiant\m e 2 copies à évaluer (préalablement anonymisées). La consigne était d'évaluer les copies et de les annoter à destination de leur auteur\m rice. Ces copies "corrigées" ont étés numérisées avant d'être rendues à leur auteur\m rice. Sur le plan du dispositif, il nous a semblé intéressant qu'un\m e étudiant\m e ait à évaluer plus d'une copie, et que chacun\m e puisse recevoir plus d'un retour sur son travail. Sur le plan de la recherche, cela permet d'obtenir de nombreuses évaluations, plusieurs pour un\m e même évaluateur\m rice et plusieurs sur un même travail.

Cependant, ces traces d'évaluation ne nous donnent qu'assez peu accès à la réflexion des étudiant\m es pendant qu'il\m elles évaluent. Nous avons pu observer que les commentaires sur les copies sont assez inégaux (en quantité et en contenu) en fonction de l'évaluateur\m rice. Ces données sont donc limitées relativement à notre problématique, et nous les avons complétées par une forme d'entretien d'explicitation. Une dizaine d'étudiant\m es volontaires ont été rassemblé\m es en binômes, que nous avons observés en train d'évaluer des copies, avec une consigne d'expliciter à voix haute ce qu'il\m elles font. Avec l'aide de l'enseignant nous avons formé des binômes d'étudiant\m es de niveaux proches, afin qu'un échange entre les deux étudiant\m es ait lieu, plutôt qu'une correction par celui ou celle de meilleur niveau. Ce dispositif a un intérêt sur le plan de la recherche : ce travail en binôme qui échange oralement et cette consigne de travail à voix haute d'explicitation de ce que font les étudiant\m es nous permettent d'accéder, au moins partiellement, à l'activité des étudiant\m es pendant qu'il\m elles évaluent leurs pairs. Ces séances de correction en binôme ont été filmées.

Précisons qu'après une rapide revue des copies, nous avons effectué une sélection pour celles que nous avons fournies aux étudiant\m es volontaires pour les entretiens. Nous avons fait en sorte que certaines de ces copies présentent des stratégies de résolution classiques et originales, et qu'elles fassent apparaître différents types d'erreurs que nous avons identifiés dans les copies.

Au final, notre corpus est composé de :
\begin{itemize}
    \item 44 copies originales,
    \item 64 copies évaluées à la maison par 30 élèves,
    \item 6 entretiens d'explicitations au cours desquels ont été évaluées entre 2 et 4 copies.
\end{itemize} 

Enfin, ce corpus a été complété d'un questionnaire (donné en annexe \ref{Qtaire}) rempli par 38 étudiant\m es, permettant de recueillir leurs impressions sur le processus d'évaluation entre pairs qu'il\m elles ont vécu. Sur le plan de la recherche, ces données nous permettent d'accéder à la perception de la tâche par les étudiant\m es et à des discours sur leur activité, possiblement sur les plans méta-cognitifs et méta-mathématiques.

\paragraph{Méthodologie d'analyse}

Après avoir sélectionné les tâches données aux étudiant\m es (section \ref{taches}) nous avons réalisé une analyse a priori des stratégies de résolution des exercices, et une analyse a priori des activités possibles lors de leur évaluation (section \ref{apriori}). Ceci nous a permis de construire une grille d'analyse des évaluations, complétée par la typologie des arguments de vérification de \cite{coppe_evaluation_2017}. Une fois les différentes données recueillies, la grille nous a permis d'identifier l'activité des étudiant\m es dans les évaluations en binômes. Ces premiers résultats nous fournissent des indications quant à l'activité des étudiant\m es en évaluation à la maison, activité à laquelle on ne peut accéder directement.

À la fin de l'expérimentation, nous avons fait passer un questionnaire aux étudiant\m es, afin de recueillir leurs impressions sur l'expérimentation et la perception qu'il\m elles ont de leur activité d'évaluation (section \ref{sectionquestionnaire}).  Nous analyserons les réponses apportées par les étudiant\m es au regard des résultats des analyses qui les précèdent.

\paragraph{Tâches proposées}
\label{taches}

Pour concevoir les exercices, nous avons d'abord sélectionné des notions familières des étudiant\m es, déjà pratiquées en cours. Les chapitres "Réels", sur le thème des bornes supérieures et inférieures, et "Ensembles" ont été choisis. Nous avons ensuite cherché des exercices qui auraient pu être posés en interrogation écrite ou oral au moment de l'expérimentation.

De plus, il était intéressant que les stratégies de résolution de ces exercices soient variées, et amènent une diversité d'erreurs potentielles dans les raisonnements. De cette manière, les étudiant\m es peuvent être confronté\m es à des raisonnements classiques et d'autres plus originaux, mais aussi à une diversité de raisonnements erronés.

Ceci a mené au choix des deux exercices suivants :

\begin{exo}
On considère l'ensemble $S=\left \{ ~ \dfrac{n-\frac{1}{n}}{n+\frac{1}{n}} ~ ~, ~ n \in \N^* \right \} $. \\

$S$ admet-il une borne supérieure ? une borne inférieure ? un minimum ? un maximum ? Si oui, les déterminer.
\end{exo}

\begin{exo}
Soient $E$ et $F$ deux ensembles, $A_1$, $A_2$ des parties de $E$, et $B_1$, $B_2$ des parties de $F$.
On note $A \times B = \{ (a,b),~a \in A \text{ et } b \in B \}$.
\begin{enumerate}
 \item Montrer que $(A_1\times B_1)\cap(A_2\times B_2)=(A_1\cap A_2)\times (B_1\cap B_2)$.
 \item A-t-on toujours $(A_1\times B_1)\cup(A_2\times B_2)=(A_1\cup A_2)\times (B_1\cup B_2)$ ?
\end{enumerate}
\end{exo}

La première question demande une preuve très formelle et requiert une bonne maîtrise de la syntaxe mathématique. La deuxième peut amener à de nombreux contre-exemples. Cependant, nous avons remarqué a posteriori que les étudiant\m es avaient très peu réussi cet exercice et qu'il avait amené à de nombreuses erreurs de résolution et de correction. L'enseignant nous a confirmé que ce type d'exercice n'était pas encore maîtrisé par la plupart des étudiant\m es.
Nous avons donc choisi de ne pas présenter en détail les analyses du second exercice.

Dans la partie suivante, nous présentons les analyses \textit{a priori} de l'expérimentation.

\subsection{Analyses a priori du problème}
\label{apriori}

Dans toute la suite, nous nous concentrons sur le premier des deux exercices :

\setcounter{exo}{0}
\begin{exo}
On considère l'ensemble $S=\left \{ ~ \dfrac{n-\frac{1}{n}}{n+\frac{1}{n}} ~ ~, ~ n \in \N^* \right \} $. \\

$S$ admet-il une borne supérieure ? une borne inférieure ? un minimum ? un maximum ? Si oui, les déterminer.
\end{exo}
Nous analysons d'abord les stratégies liées à la tâche de résolution du l'exercice, puis nous analysons la tâche d'évaluation associée.

\subsubsection{Analyse a priori des stratégies de résolution}

Nous présentons ici les principales stratégies de résolution du premier exercice. Il peut être résolu en adoptant un point de vue "suites numériques" ou un point de vue "ensembles" sur les objets en jeu. Bien entendu, des méthodes hybrides entre ces deux points de vues sont envisageables.

Notons pour tout $n \in \N^*$, $s_n = \dfrac{n-\frac{1}{n}}{n+\frac{1}{n}}$.

\paragraph{En adoptant le point de vue des suites numériques,} on montre que la suite $(s_n)$ est croissante. Ceci est possible par exemple :
\begin{itemize}
    \item en étudiant la différence $s_{n+1}-s_n$;
    \item en étudiant le quotient $s_{n+1}/s_n$;
    \item en étudiant la croissance de la fonction $f : x \in \R \mapsto \dfrac{x-\frac{1}{x}}{x+\frac{1}{x}}$.
\end{itemize}

L'ensemble des termes de la suite est donc minoré par le premier terme $s_0$, ce qui prouve l'existence de la borne inférieure et du minimum de l'ensemble $S$, et permet de trouver sa valeur : $\inf(S) = s_0=0$. Pour montrer l'existence de la borne supérieure de $S$, on peut commencer par montrer que $(s_n)$ converge vers une limite finie : 

\begin{equation*}
    \text{Pour tout }n \in \N^*, \text{ on a : } s_n = \dfrac{n-\frac{1}{n}}{n+\frac{1}{n}} = \dfrac{1-\frac{1}{n^2}}{1+\frac{1}{n^2}} \text{ donc }s_n \underset{n \to + \infty}{\longrightarrow} 1.
\end{equation*}

Le théorème de la limite monotone permet de conclure que l'ensemble $S$ a une borne supérieure, égale à 1. Il reste alors à montrer que cette borne supérieure n'est pas un maximum. Supposons par l'absurde qu'il existe un entier $n \in \N^*$ tel que $\dfrac{n-\frac{1}{n}}{n+\frac{1}{n}} = 1$. Cette égalité est équivalente à $n-\frac{1}{n} = n+\frac{1}{n}$, elle même équivalente à l'égalité $-1=1$, ce qui est une contradiction. Aucun élément de l'ensemble $S$ n'est donc égal à 1, la borne supérieure n'est pas un maximum.

\paragraph{Retrouvons ce résultat avec un point de vue "ensembles".} On peut encadrer directement les éléments de l'ensemble. En effet, pour tout $n \in \N^*$, on a :

\begin{equation*}
    n-\frac{1}{n} \geq 0 \text{ et } n+\frac{1}{n} >0 \text{ donc } \dfrac{n-\frac{1}{n}}{n+\frac{1}{n}} \geq 0
\end{equation*}

et de plus

\begin{equation*}
    n-\frac{1}{n} < n+\frac{1}{n} \text{ donc } \dfrac{n-\frac{1}{n}}{n+\frac{1}{n}} < 1.
\end{equation*}

Ceci montre que pour tout $n \in \N^*$, $0 \leq s_n < 1$. L'ensemble $S$ admet donc une borne inférieure et une borne supérieure, qui vérifient $\inf(S) \geq 0$ et $\sup(S) \leq 1$. Or $s_0 = 0$, donc $\inf(S)=\min(S)=0$. Pour montrer $\sup(S) = 1$, on peut utiliser la caractérisation séquentielle de la borne supérieure, et montrer comme ci-dessus que $s_n \underset{n \to +\infty}{\longrightarrow} 1$. On peut également faire appel à la caractérisation "avec des epsilons" de la borne supérieure. Soit donc $\varepsilon >0$, montrons qu'il existe $n_0 \in \N^*$ tel que $1-\dfrac{n_0-\frac{1}{n_0}}{n_0+\frac{1}{n_0}} < \varepsilon$. Or si $n \in \N^*$, on a les équivalences :

\begin{equation*}
    1-\dfrac{n-\frac{1}{n}}{n+\frac{1}{n}}  < \varepsilon \quad \Longleftrightarrow \quad  \frac{2}{n^2+1}  < \varepsilon \quad  \Longleftrightarrow  \quad \frac{2}{\varepsilon} -1 < n^2.
\end{equation*}

Si $\frac{2}{\varepsilon} -1 \leq 0 $, $n_0=1$ convient. Sinon, on pose $n_0$ la partie entière supérieure de $\sqrt{\frac{2}{\varepsilon} -1}$, qui vérifie bien l'inégalité demandée. On a donc montré que la borne supérieure de $S$ est égale à 1. De plus, l'inégalité $s_n <1$ vraie pour tout $n \in \N^*$ permet de conclure que l'ensemble $S$ n'a pas de maximum. 

Notons qu'on peut réécrire les éléments de l'ensemble grâce à l'égalité $s_n = 1 - \frac{2}{n^2+1}$. Ceci permet de montrer rapidement que la suite $(s_n)$ est croissante, qu'elle tend vers 1, ou encore que 1 n'est pas un élément de l'ensemble.

\paragraph{Retour sur nos contraintes} Cet exercice vérifie donc les contraintes que nous nous sommes données : les visions "suites" ou "ensembles" permettent de varier les types de raisonnements rencontrés dans les copies. Au sein de la vision "suites", différentes démonstrations sont possibles pour montrer une même propriété, comme la croissance de la suite. Nous nous attendons donc à rencontrer divers raisonnements et permettant des activités différentes lors des évaluations.

\subsubsection{Analyse a priori des évaluations}

Une analyse \textit{a priori} de l'activité des étudiant\m es et des arguments d'évaluation est présentée dans cette partie, et a servi de grille d'analyse des évaluations. Cette analyse s'articule en plusieurs niveaux. Cinq grands types de d'activités sont définis, chacun étant décliné en deux sous-niveaux, qui précisent et détaillent les niveaux précédents. Une description des arguments d'évaluation adaptée de \cite{coppe_evaluation_2017} est ensuite utilisée.

\paragraph{Classification générale des évaluations}
 
Nous avons distingué trois catégories d'évaluations : celles portant sur le fond du raisonnement mathématique, celles sur la forme, et celles ne portant pas sur le raisonnement mathématique. Au sein de ces catégories, on retrouve les types d'activités suivantes :

\begin{enumerate}[label=\textbf{Catégorie \arabic*}, left=0.5cm]
    \item Activités portant sur l'évaluation d'une preuve proposée dans la copie ;
    \item Activités portant sur l'évaluation d'une preuve envisagée par l'étudiant\m e corrigé\m e, selon l'évaluateur\m rice ;
    \item Activités pourtant sur l'évaluation du résultat d'une preuve proposée dans la copie ;
    \item Activités pourtant sur l'évaluation de la rédaction ou de la syntaxe d'une preuve ;
    \item Activités ne portant pas sur l'évaluation du raisonnement mathématique.
\end{enumerate}

Notons qu'une distinction est faite entre l'évaluation d'une preuve telle qu'elle est proposée dans la copie, et l'évaluation d'une interprétation de cette preuve, faite par l'évaluateur\m rice.



\paragraph{Détails du premier niveau de classification}

Chaque catégorie a été déclinée en plusieurs familles d'activités. Par exemple, lorsque les étudiant\m es évaluent une preuve proposée dans la copie - catégorie 1 -, il\m elles sont amené\m es à :

\begin{enumerate}[label=\Alph*]
\addtocounter{enumi}{+25}
    \item Lire la preuve;
\end{enumerate}
\begin{enumerate}[label=\Alph*]
    \item Essayer de comprendre la preuve ;
    \item Évaluer la validité de la preuve ;
    \item Évaluer la vraisemblance de la preuve ;
    \item Corriger la preuve.
\end{enumerate}

Ces activités peuvent apparaître simultanément ou séparément. Par exemple, certains étudiant\m es lisent toute la copie, puis évaluent la validité des raisonnements tout en les corrigeant. D'autres pourront lire et évaluer la validité dans un même temps, puis corriger chaque exercice.

\paragraph{Détails du deuxième niveau de classification}

Au sein de chaque famille, nous avons enfin identifié plusieurs activités que les étudiant\m es peuvent mettre en place lors de l'évaluation. Détaillons l'ensemble \textit{\'Evaluer la validité de la preuve}. Les étudiant\m es vont être amené\m es à valider ou invalider la preuve complète ou seulement une partie de celle-ci (un calcul, une implication \dots). Les étudiant\m es s'interrogent aussi sur la démarche, "l'idée de la preuve", entamée dans la copie. Nous incluons dans cette famille d'activités les questionnements sur la validité de la preuve, y compris lorsqu'il n'aboutissent pas.

Décrivons à présent la famille \textit{Essayer de comprendre la preuve}. Les étudiant\m es peuvent chercher à identifier les différentes étapes de la preuve qu'il\m elles évaluent. Dans notre exercice, un\m e étudiant\m e pourra d'abord montrer que la borne supérieur existe en majorant les éléments de l'ensemble, puis montrer que la suite est croissante pour conclure sur la valeur de la borne supérieure. Pendant ce processus, les évaluateur\m rices sont amené\m es à identifier les concepts mathématiques mobilisés, implicitement ou explicitement, dans la copie. Lors des entretiens, les étudiant\m es peuvent expliquer le raisonnement mené dans la copie (ou une partie de celui-ci) à leur camarade. Enfin, certain\m es étudiant\m es peuvent anticiper la suite du raisonnement en même temps qu'il\m elles lisent et détaillent la preuve.

\paragraph{Une classification des arguments de vérification}

Pour enrichir l'analyse de l'activité des étudiant\m es, nous avons introduit un deuxième niveau de description, basé sur la classification de \citep{coppe_evaluation_2017}. Dans les entretiens en binômes, nous identifions dans le discours des étudiant\m es des arguments à visée de validation ou d'invalidation de preuves, de correction ou encore d'explication d'une preuve. Pour formuler ou étayer ces arguments, les étudiant\m es peuvent être amené\m es à faire appel à des références comme le cours, la feuille de correction, ou leur propre expérience.

Nous appelons argument ou référence \textit{externes} ceux ne faisant pas intervenir de savoir ou de savoir faire mathématique. On a a priori listé :

\begin{enumerate}
    \item les références à l'énoncé de l'exercice;
    \item les références à une expérience personnelle de l'étudiant\m e;
    \item les références à la mémoire du savoir enseigné;
    \item les références aux éléments de corrections.
\end{enumerate}

Les arguments ou références nécessitant un raisonnement mathématique sont dits \textit{internes}.
On peut identifier :

\begin{enumerate}
    \item les vérifications d'éléments techniques :
    \begin{itemize}
        \item vérifications d'un calcul à la main;
        \item vérifications d'une preuve à la main;
        \item vérifications de la présence (ou non) des hypothèses d'un théorème;
        \item production d'exemples ou de contre-exemples pour mettre à l'épreuve un énoncé;
    \end{itemize}
    \item les références aux résultats et méthodes du cours ou d'une source experte.
\end{enumerate}

De manière similaire à l'article de référence, nous avons pensé aux vérifications par changement de cadre, mais nous ne les avons jamais observées en pratique.

\subsubsection{Retour sur les choix expérimentaux}

Rappelons que nous avons choisi, suite à l'expérimentation, de nous concentrer sur les évaluations du premier exercice.
L'enseignant de la classe a confirmé que les étudiant\m es avaient à disposition les connaissances pour résoudre l'exercice. L'analyse a priori de cet exercice montre la diversité des résolutions possibles, faisant appel aux domaines des suites numériques et des ensembles. Les étudiant\m es peuvent donc mobiliser différentes notions, ce qui pourra produire des résolutions puis des évaluations diverses.

L'analyse a priori de l'évaluation nous a permis de construire une grille de codage, utilisée pour analyser le travail réalisé par les binômes. Sa version finale est présentée en annexe \ref{annexegrillenalayse}. Cette grille est adaptée aux entretiens, durant lesquels les étudiant\m es verbalisent leurs réflexions au cours de l'évaluation. En revanche, il parait compliqué d'identifier toutes ces activités uniquement à partir des copies annotées. Cette grille est donc amenée à être adaptée et allégée pour être pertinente dans ce cadre. Par conséquent, dans la suite de ce texte, nous analyserons d'abord les entretiens, et nous chercherons ensuite des traces, dans les copies annotées, des activités que nous aurons identifiées préalablement.

\section{Résultats des analyses de l'activité des étudiant\m es}

Avant de présenter les analyses des entretiens et des copies, nous donnons en première partie quelques éléments de repère concernant la résolution et l'évaluation des exercices par les étudiant\m es.

\subsection{Premières remarques sur la résolution et l'évaluation des exercices}

\paragraph{Sur la résolution des exercices}

Avant l'évaluation par les étudiant\m es, nous avons consulté l'ensemble des copies afin de sélectionner celles que nous avons donné à corriger en entretiens. L'analyse de ces copies nous permet de faire un état des lieux des stratégies utilisées par les étudiant\m es et des erreurs les plus courantes.

La stratégie la plus fréquemment employée (34 étudiant\m es) passe par l'étude de la croissance de la suite. Pour cette stratégie, les étudiant\m es sont réparti\m es en :

\begin{itemize}
    \item calcul de la dérivée de la fonction associée : 25 étudiant\m es;
    \item étude de $s_{n+1}-s_n$ : 7 étudiant\m es;
    \item étude de $s_{n+1}/s_n$ : 1 étudiant\m e;
    \item décomposition en éléments simples du terme général : 1 étudiant\m e (qui étudie ensuite $s_{n+1}-s_n$).
\end{itemize}

Parmi les étudiant\m es qui n'ont pas utilisé cette stratégie, une partie argumente que le premier terme de la suite est la borne inférieure de l'ensemble et que la limite de la suite est la borne supérieure, sans plus de justifications. Deux étudiant\m es tentent d'utiliser la caractérisation séquentielle de la borne supérieure, et un étudiant\m e cite la caractérisation "avec des $\varepsilon$", sans l'utiliser ensuite. Au final, aucun\m e de ces étudiant\m es n'a produit de preuve valide que la borne supérieure de l'ensemble est 1. Enfin, certain\m es étudiant\m es minorent la suite par 0 et montrent l'égalité $s_1=0$ pour conclure quant à la valeur de la borne inférieure.

Le premier examen des copies nous a également permis d'identifier les erreurs les plus fréquentes (présentes dans plus de 5 copies) dans le premier exercice. Il s'agit de :

\begin{itemize}
    \item si $s_n \to 1$ lorsque $n \to + \infty$; alors $(s_n)$ est majorée par 1;
    \item si $s_n \to 1$ lorsque $n \to + \infty$, alors $\sup(S)=1$;
    \item si $s_1 = 0$, alors 0 minore $S$;
    \item dérivation de la fonction $f : \N^* \to \R$ associée à la suite $(s_n)$.
\end{itemize}

\paragraph{Sur l'évaluation des copies}

\begin{figure}
    \centering
    \includegraphics[width=0.8\textwidth]{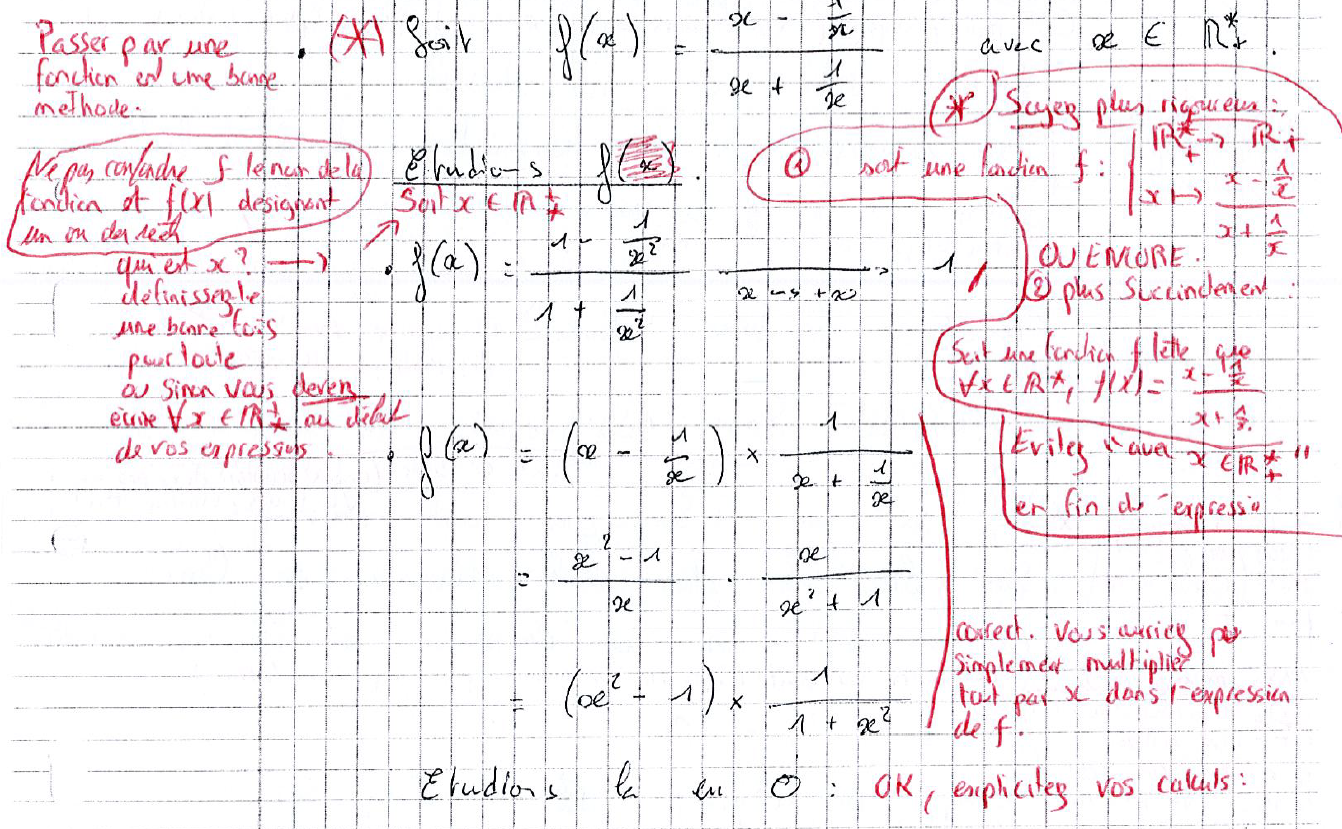}
    \caption{Un exemple de correction de l'exercice 1 très détaillée}
    \label{cop44}
\end{figure}

\begin{figure}
    \centering
    \includegraphics[width=0.8\textwidth]{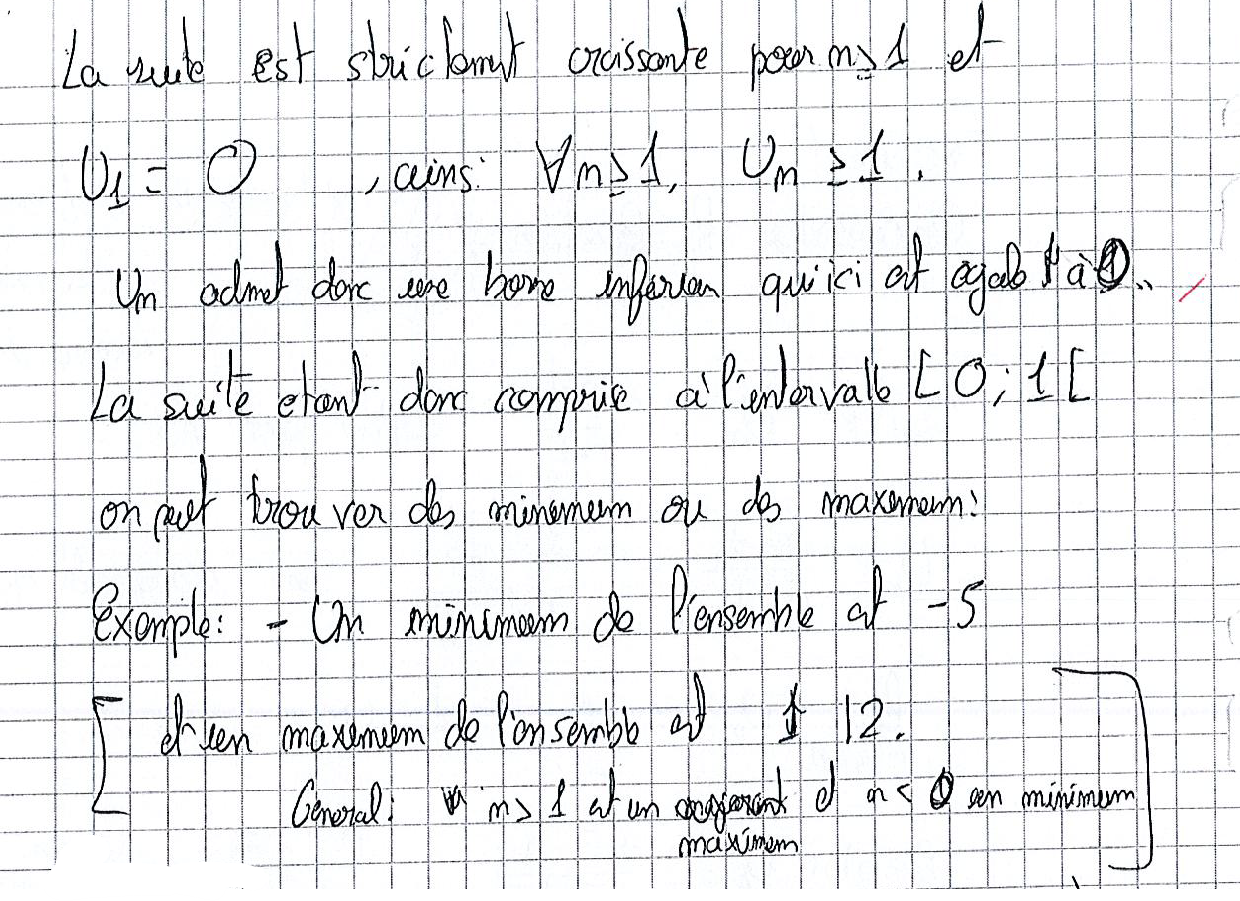}
    \caption{Un exemple de correction (incorrecte) de l'exercice 1 avec très peu de détails}
    \label{cop18}
\end{figure}

30 étudiant\m es ont chacun\m e évalué 2 copies en travail à la maison. À l'examen de ces évaluations, on observe une grande variabilité en fonction de l'évaluateur\m rice, tant au niveau de la validité de la correction que dans les annotations. Certain\m es étudiant\m es donnent de nombreux commentaires et conseils, tandis que d'autres se contentent de valider ou d'invalider les résultats. Des exemples de ces comportements sont présentés dans les figures \ref{cop44} et \ref{cop18}. À la fin de l'expérimentation, dans une perspective formative, nous avons donc examiné toutes les évaluations afin de les compléter et les corriger si nécessaire.

Concentrons nous à présent sur les résultats des entretiens. Rappelons que lors de ces entretiens, des binômes d'étudiant\m es ont évalué ensemble entre 2 et 4 copies, avec la consigne d'expliciter oralement leurs discussions, réflexions et commentaires sur l'évaluation des copies.

\subsection{Analyse des entretiens}

Pour l'analyse, nous avons travaillé à partir de transcriptions des entretiens. Nous avons d'abord codé l'activité des étudiant\m es de chaque binôme à l'aide de la grille d'analyse présentée ci-dessus.
Sur cette base, nous réalisé une synthèse de l'activité de chaque binôme, puis essayé de repérer les aspects les plus caractéristiques de leurs activités. Nous présentons ici une synthèse de nos résultats.

Lors de l'analyse des entretiens, nous avons dégagé sept axes caractéristiques de l'activité des étudiant\m es. Il s'agit de la compréhension de la preuve, l'évaluation de sa validité ou de sa vraisemblance, l'attention portée au résultat de la preuve, l'attention portée à la démarche de résolution empruntée dans la copie, les corrections apportées par les évaluateur\m rices, les remarques concernant la syntaxe et la rédaction et leurs liens avec l'invalidation de la preuve, et enfin les jugements de valeur. Pour illustrer nos résultats, nous choisissons de présenter dans la suite plus en détails trois de ces axes, qui nous semblent au coeur des enjeux de l'entrée dans les mathématiques du supérieur :  compréhension de la preuve, validation et invalidation, corrections apportées. Nous reviendrons ensuite sur les arguments et ressources mobilisés par les étudiant\m es.

Pour les exemples que nous donnons dans cette section, nous notons les binômes \textbf{B1}, \textbf{B2}, \textbf{B3}, \textbf{B4}, \textbf{B5}, \textbf{B6}.

\paragraph{Compréhension de la preuve}

Lors de l'évaluation d'une copie, les étudiant\m es observé\m es commencent par lire et essayer de comprendre les preuves proposées. Les activités de lecture et de validation d'une preuve s'effectuent la plupart du temps en parallèle. Il arrive qu'une phase uniquement d'explicitation et d'explication de la preuve précède toute validation ou invalidation, comme pour le binôme \textbf{B4}.

\begin{quote}
    \og Je sais pas si c'est juste, mais c'est ça qu'il fait. \fg
\end{quote}

Lors de l'examen de la preuve, les étudiant\m es identifient et détaillent les différentes étapes qui la compose, comme dans la citation suivante du binôme \textbf{B4}.

\begin{quote}
    \og Après il a démontré qu'elle était croissante, du coup son minimum c’est sa première valeur vu qu'elle est strictement croissante. \fg
\end{quote}


Dans cette phase, les étudiant\m es sont amené\m es à identifier précisément les concepts mathématiques utilisés, qui parfois ne sont pas explicités dans la copie. Nous avons régulièrement observé cette situation dans les entretiens, lorsque les étudiant\m es évalué\m es utilisent le théorème de la limite monotone sans le citer, ou énoncent un théorème proche. Il arrive alors que les évaluateur\m rices interprètent le raisonnement mené dans la copie à la lueur des théorèmes qu'il\m elles connaissent.

\begin{quote}
    \og (\textbf{B3}) En fait il a utilisé la limite monotone, mais dans l'autre sens. \fg
\end{quote}



Ceci est aussi une occasion pour les étudiant\m es de se rappeler précisément des définitions ou théorèmes du cours, afin de comprendre si les concepts mobilisés le sont de la bonne manière. Par exemple, le binôme \textbf{B6} débat à deux reprises sur les différentes caractérisations de la borne supérieure vue en cours : définition, caractérisation séquentielle et caractérisation "avec des epsilons". Notons également que contrairement aux autres binômes, les étudiant\m es du binôme \textbf{B1} tentent peu d'explicitations et explications détaillées des raisonnements à l'oral. Il\m elles se contentent le plupart du temps de lire la preuve à voix haute, et notent quand les démarches mises en \oe uvre dans la copie sont similaires aux leurs. On remarquera  dans la suite que ce binôme utilise beaucoup de références externes pour les validations et les corrections des exercices.

Par ailleurs, s'il\m elles n'ont pas réussi à comprendre un raisonnement ou un résultat annoncés dans la copie, les évaluateur\m rices peuvent interpréter \og ce que l'étudiant\m e évalué\m e a voulu faire dans la preuve \fg. Dans une copie corrigée par le binôme \textbf{B5}, il est écrit "$u(x) \in [0,1[$", où $u$ est la fonction associée à la suite introduite dans l'exercice 1. L'étudiant\m e justifie correctement l'encadrement $u(x) \in [0,1]$, mais sans prouver l'ouverture en 1, et sans faire le lien avec la borne supérieure de la suite. Le binôme a interprété que l'ouverture de l'intervalle a été citée pour justifier que 1 n'est pas le maximum de l'ensemble. Ces comportements font donc partie de la catégorie \textit{Activités portant sur l'évaluation d'une preuve envisagée par l'élève corrigé, selon le correcteur}.

On remarque enfin que les étudiant\m es du binôme \textbf{B2} ont la spécificité de souvent anticiper les raisonnements à venir dans la copie. En conséquence, si la copie emprunte le chemin de résolution envisagé par les étudiant\m es, la preuve est validée même si elle est incorrecte. Réciproquement, une preuve correcte mais qui utilise une autre méthode peut être invalidée. Ces biais influent donc sur la validité des évaluations produites par ce binôme. On trouve ci-dessous deux exemples de ces anticipations. 

\begin{quote}
    \textit{Face à une copie dans laquelle il est seulement mis "$u_n$ croissante et converge vers 1".}
    
    \og Et du coup si elle converge vers 1, ça veut dire qu'en même temps 1 c'est le max et le sup ça voudrait dire. Est ce qu'il l'a bien expliqué ? \fg
\end{quote}

\begin{quote}
    \og Et est ce qu'il dit que 1 appartient à l'ensemble ? Ça vient peut être après. \fg
\end{quote}

En conclusion, on observe toujours un essai de compréhension des preuves proposées dans la copie évaluée, qui est souvent parallèle à une première validation. Les étudiant\m es lisent la preuve et identifient explicitement les étapes qui la constituent. Certains relèvent aussi les concepts et énoncés mathématiques en jeu.

\paragraph{Validation et invalidation, vraisemblance}

Au cours de l'évaluation, tous\m tes les étudiant\m es discutent de la validité des preuves, résultats ou démarches qu'il\m elles rencontrent dans les copies. Nous présentons ici les caractéristiques de l'évaluation de la validité des preuves. Différentes situations peuvent amener les étudiant\m es à invalider un élément d'une copie.

\begin{figure}
    \centering
    \includegraphics[width=0.8\textwidth]{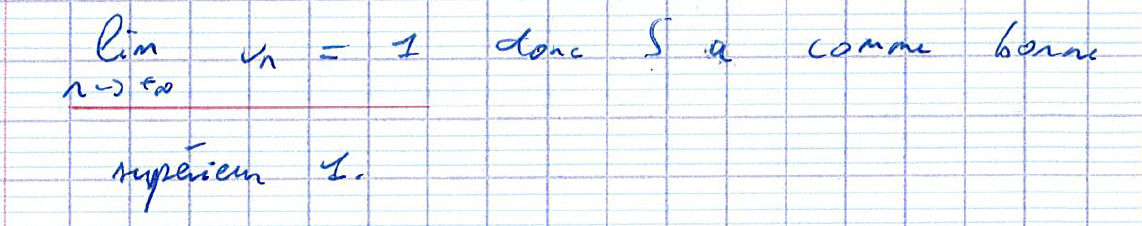}
    \caption{Utilisation implicite de l'énoncé "toute suite convergente converge vers la borne supérieure de l'ensemble de ses éléments"}
    \label{cop34}
\end{figure}

Premier cas : Lorsque, implicitement ou explicitement, une propriété avec quantification universelle invalide est utilisée dans la copie, les évaluateur\m rices produisent des contre-exemples. La plupart du temps, cela se produit quand l'étudiant\m e corrigé\m e utilise un théorème du cours sans en vérifier toutes les hypothèses. Par exemple, dans la figure \ref{cop34}, l'étudiant\m e semble utiliser l'énoncé "toute suite convergente converge vers la borne supérieure de l'ensemble de ses éléments".

Deux autres raisonnements invalides ont été relevés à plusieurs reprises par les binômes : l'utilisation de la propriété "toute partie non vide et majorée de $\R$ admet une borne supérieure" sans expliciter que l'ensemble est non vide; et l'utilisation du théorème de la limite monotone sans montrer que la suite est croissante ou qu'elle est majorée. Nous donnons ci-dessous une situation où le binôme \textbf{B3} donne un contre-exemple pour invalider la phrase "Pour $n=1, s_1 = 0$ est la borne inférieure de $S$".

\begin{quote}
    \og Si elle était décroissante elle serait pas minorée par son premier terme.
\fg \end{quote}





Deuxième cas : les étudiant\m es relèvent une erreur lorsqu'un résultat leur paraît trop peu justifié. Nous avons fréquemment observé cette situation dans le premier exercice, où des étudiant\m es citent que "la limite 1 n'est pas atteinte, donc 1 n'est pas le maximum de l'ensemble" sans démonstration. Cette affirmation est alors invalidée par les évaluateur\m rices la plupart du temps. Lorsqu'une preuve est peu justifiée, les étudiant\m es sont amené\m es à discuter de si elle est "assez justifiée" pour qu'ils la valident. Le binôme \textbf{B4} dit par exemple :

\begin{quote}
    \og --- La première partie de l'exercice 1 elle est bien écrite, mais la deuxième partie sur le max et la borne sup c'est un peu rapide, ça aurait pu être un peu plus détaillé.

    --- Mais c'est la bonne réponse.

    --- Par exemple là "par somme puis quotient", je trouve pas ça assez détaillé.

    --- Après bon là il a pas donné le nom du théorème, mais il a dit qu'elle tendait vers sa borne sup, il manque le nom. Il faut dire pourquoi il dit ça.
\fg \end{quote}

Les évaluateur\m rices se posent donc la question des justifications nécessaires pour qu'une preuve soit considérée comme valide. Ils se demandent si un théorème non cité mais utilisé correctement est valide, ou quels types de justifications sont nécessaires quand on mène des calculs. 




Troisième cas : les étudiant\m es discutent de la validité de propriétés des objets spécifiques à l'exercice. Par exemple, dans une copie corrigée par le binôme \textbf{B4}, il est affirmé que la fonction associée à la suite $s_n$ est parfois décroissante, ce qui est invalidé par le binôme évaluateur.

On observe donc plusieurs motifs d'invalidation au cours des entretiens. Elle peut porter sur les propriétés d'un objet utilisé dans le raisonnement, ou sur un théorème annoncé. Les évaluateur\m rices peuvent alors justifier l'invalidation d'un énoncé en trouvant un contre-exemple. Les étudiant\m es sont aussi amené\m es à s'interroger sur les justifications nécessaires ou suffisantes pour valider une preuve ou un calcul. Notons aussi que certaines preuves sont évaluées comme invalides soit par rejet de la démarche entamée, soit parce que le résultat obtenu est incorrect.

Du côté de la validation d'une preuve, elle est le plus souvent directe et sans justification particulière : \og Ça, c'est juste. \fg. La validation est parfois faite au regard des éléments de corrections qui ont été fournis aux étudiant\m es. Pour les éléments techniques comme les calculs ou les tableaux de signe, la vérification à la main est systématique pour le binôme \textbf{B2}, et présente ponctuellement pour les binômes \textbf{B3}, \textbf{B5} et \textbf{B6}. Il arrive aussi que les étudiant\m es ne vérifient pas les calculs, mais les valident car il\m elles se souviennent avoir obtenu le même résultat lors de leur résolution de l'exercice. Comme nous l'avons remarqué précédemment, les étudiant\m es du binôme \textbf{B1} font souvent référence à leur propre résolution des exercices pour valider les preuves qu'il\m elles évaluent :

\begin{quote} \og 
    En gros c'est ce qu'il a fait, il a calculé le terme 1, ça fait 0. Moi après j'ai fait comme ça.
\fg \end{quote}

On observe donc que les validations sont moins argumentées que les invalidations, et font parfois référence aux expériences personnelles des étudiant\m es et aux éléments de correction qui leur ont été fournis.

\paragraph{Corrections apportées}

Après avoir repéré une erreur dans la copie qu'ils évaluent, les binômes observés proposent des éléments de correction, sous forme d'annotations.

Notons d'abord que la correction des erreurs de définition est systématique. La correction des calculs n'est quant à elle pas toujours présente. Certains binômes décortiquent les calculs pour chercher les éventuelles erreurs et les corriger, d'autres se contentent de noter quand le résultat est faux.

Nous avons identifié différents types de correction. Les étudiant\m es peuvent compléter le raisonnement entrepris dans la copie, en ajoutant les arguments ou les justifications pour qu'il devienne valide. Cela amène les étudiant\m es à réfléchir aux arguments nécessaires et aux arguments suffisants à donner en correction. Dans le dialogue suivant du binôme \textbf{B3}, un\m e étudiant\m e propose une certaine méthode pour montrer que la suite est minorée. L'autre étudiant\m e corrige de manière plus générale, en affirmant qu'il faut montrer qu'elle est minorée par 0. Les étudiant\m es discutent ensuite de s'il est nécessaire d'étudier la monotonie de la suite pour montrer que la borne inférieure est 0.

\begin{quote} \og 
    --- Il faudrait dire que si jamais cette suite elle est croissante, et que du coup elle est minorée par son premier terme. Du coup ça c'est une partie de ça.
    
    --- Oui il faut d'abord qu'il montre qu'elle est minorée par 0. Je mets "il manque 0 minore S".

    --- Oui et puis aussi l'étude de la suite quoi. Savoir si elle est croissante ou décroissante.
    
    --- Ben je veux dire pour ça, pour cette partie en particulier, il a pas, enfin ... pour montrer que 0 minore S tu utilises l'étude de la suite, non ?
\fg \end{quote}

La correction du binôme \textbf {B4} ci-dessous nous donne un second exemple. Les étudiant\m es commencent par invalider l'argument fourni en proposant un contre-exemple : une suite stationnaire à 1 vérifierait les mêmes conditions de croissance, mais aurait un maximum. En regardant les éléments de correction, il\m elles affirment ensuite qu'il faudrait encadrer la suite pour montrer que la borne supérieure n'est pas un maximum. La méthode proposée est donc suffisante. Cependant, pour conclure, il est seulement nécessaire de montrer que la borne supérieure n'est pas incluse dans l'ensemble.

\begin{quote} \og 
    --- Il a mis limite de f(x) = 1. Du coup ça a pas de max, mais c'est pas vraiment un max quand même.

    --- Oui c'est vrai qu'elle peut être stationnaire.
    
    --- Il fallait qu'il encadre la suite.
\fg \end{quote}

Ensuite, les binômes peuvent corriger une preuve en réorganisant les arguments, pour "remettre la preuve dans le bon ordre" comme dans l'exemple ci-dessous. La compréhension de l'ordre et de la bonne structure d'une preuve mathématique peut être une difficulté pour les étudiant\m es en début de l'enseignement supérieur, et on retrouve en effet plusieurs erreurs de ce type dans les copies.

\begin{quote}
    (\textbf{B5}) \og Genre dès le début il aurait pu commencer avec la limite, il aurait pu commencer avec le min et la borne inf, ensuite embrayer avec la croissance. Parce que ce qu'il a fait c'est pas mauvais, il a donné des arguments mais il a vendu la peau de l'ours avant de l'avoir tué en gros.
\fg \end{quote}

Enfin, les évaluateur\m rices peuvent proposer une autre piste de résolution que celle envisagée dans la copie. Il arrive plus rarement que des étudiant\m es proposent une autre méthode de résolution qu'ils jugent plus simple, après avoir validé celle présentée dans la copie. Dans l'exemple suivant, le binôme \textbf{B2} pense que le théorème invoqué dans la copie ne permet pas de résoudre l'exercice comme annoncé. Plutôt que de fournir une correction se basant sur la caractérisation séquentielle de la borne supérieure comme entamé dans la copie, les étudiant\m es affirment qu'il faut utiliser le théorème de la limite monotone pour conclure. 

\begin{quote} \og
    Sup de S par caractérisation séquentielle ? Théorème des limites monotones ? Ouais c'est juste ... C'est le bon résultat avec une enfin ... Non non c'est juste le TLM \textit{(Note : théorème de la limite monotone)}.
\fg \end{quote}


Les corrections concernent donc plusieurs contenus mathématiques : définitions, calculs, preuves \dots Les étudiant\m es peuvent compléter une preuve en rajoutant des éléments et des justifications pour la rendre valide. Notons que certains binômes cherchent à donner des arguments suffisants minimaux dans la correction, et on peut supposer que ces étudiant\m es font preuve d'une compréhension avancée du raisonnement déductif et des concepts et énoncés mathématiques en jeu. D'autres, en revanche, complètent la preuve avec une condition suffisante non-nécessaire, comme un binôme qui propose de montrer que la suite est croissante pour prouver qu'elle est minorée. Enfin, les étudiant\m es ont régulièrement tenté de remettre dans le bon ordre les arguments identifiés dans la copie, pour reconstruire une preuve valide.

Nous nous concentrons maintenant sur les arguments et ressources mobilisés par les étudiant\m es pendant les différentes phases de l'évaluation.

\paragraph{Arguments internes et externes}

Au cours de l'évaluation des copies, les étudiant\m es font appel à certains arguments, à visée de validation de preuves ou de commentaire des copies. Rappelons qu'un argument ou une référence est dit \textit{externe} s'il ne fait pas intervenir de savoir ou savoir-faire mathématique. Ceux qui requièrent un raisonnement mathématique sont dits \textit{internes}.

L'utilisation de contre-exemples est présente pour tous les binômes, principalement pour invalider des arguments ou des résultats. Ces réflexions (arguments internes) peuvent amener à une compréhension plus profonde des théorèmes de leur cours. En effet, elles permettent de comprendre la raison d'être de certaines hypothèses des théorèmes que les étudiant\m es utilisent dans leur pratique mathématique. Cependant, certain\m es étudiant\m es utilisent peu de contre-exemples pour invalider des preuves, au profit d'arguments et références externes.

Du côté des arguments externes, tous les binômes font des références à des expériences personnelles. La plupart du temps, elles visent à commenter la preuve qu'ils évaluent, et à la comparer avec leur propre résolution de l'exercice. Par exemple, les étudiant\m es du binôme \textbf{B5} discutent de leur résolution de l'exercice quand il\m elles essaient de comprendre la preuve engagée dans la copie. Il\m elles disent alors :

\begin{quote} \og 
    Par contre là il a fait mieux que moi. Moi j'avais pas mis le "et", j'avais mélangé ensemble et ... du coup en lisant sa copie j'ai appris de mon erreur aussi. 
\fg \end{quote}

On voit déjà ici que l'évaluation de certaines copies, dont la stratégie de résolution est proche de celle engagée par les évaluateurs, peut leur permettre de prendre du recul sur leur propre pratique. Notons que les étudiant\m es du premier binôme font de nombreuses référence de ce type pour valider des preuves, sans plus d'argumentation. Ces étudiant\m es n'utilisent que peu d'arguments internes.

Il nous a semblé que lors de la validation, l'invalidation ou la correction d'une preuve, les étudiant\m es utilisent des arguments externes lorsqu'il\m elles ne sont pas à l'aise avec les concepts en jeu, et ont des difficultés à mobiliser des arguments internes. Au niveau des mathématiques du supérieur, il paraît souhaitable que les étudiant\m es apprennent à mobiliser des arguments internes pour évaluer les preuves, arguments qui peuvent aussi servir lors de la construction de nouvelles preuves.

\paragraph{Ressources utilisées}

On observe dans ce recueil de données que les étudiant\m es font aussi appel à des ressources extérieures lors de l'évaluation de copies de pairs. L'utilisation de ces ressources leur permet de valider, invalider, ou commenter les copies qu'il\m elles évaluent. On trouve naturellement quelques références aux cours auxquels les étudiant\m es ont assisté sur une notion, mais aussi aux supports de correction qui ont été fourni durant l'expérimentation.

On note d'abord que les étudiant\m es font peu de références explicite à leur cours, et lorsqu'il\m elles le font, les propriétés mentionnées sont parfois incorrectes. Par exemple, lors de la correction de l'exercice 2, le binôme \textbf{B1} affirme que le produit cartésien se développe sur l'union (ce qui est à la fois problématique car c'est une propriété fausse, et parce que cette propriété est en fait l'objet même de la question). Seul\m e un\m e étudiant\m e, dans le binôme \textbf{B2}, semble faire référence de manière précise à son cours. Les références au cours sont utilisées pour valider des preuves, mais aussi pour vérifier des définitions énoncées dans les copies. La situation semble donc différente de celle de recherche de problèmes, dans laquelle on s'attend à ce que le cours soit un support à la réflexion. Notons cependant que les étudiant\m es ont déjà cherché l'exercice, et qu'il\m elles ont pu à cette occasion se remémorer le cours sur le sujet. Nous avons aussi remarqué que la vérification de la présence des hypothèses d'un théorème n'est pas systématique chez les étudiant\m es qui évaluent, ou n'est pas toujours aboutie. On peut supposer que cela est dû à une connaissance approximative des énoncés du cours, ou à un rapport aux théorèmes où l'importance des hypothèses est négligée. Cela peut amener les étudiant\m es à valider des preuves qui sont incorrectes. Le cours ne semble donc pas être un support utilisé de manière explicite et précise.

Ensuite, quelques références à la feuille des éléments de correction sont faites. Les étudiant\m es l'utilisent principalement pour valider une démarche, si celle-ci apparaît dans la feuille. Pour certain\m es étudiant\m es, cela est utile dans la phase de compréhension de la preuve qu'il\m elles évaluent : ils s'en servent pour identifier précisément les méthodes engagées dans la copie.

Les étudiant\m es mobilisent donc peu les ressources extérieures à leur disposition dans l'évaluation de leur copie.

Nous avons détaillé dans cette section trois principaux types d'activités observés dans le processus d'évaluation par les pairs : compréhension d'une preuve, validation et invalidation, et correction. Nous avons aussi remarqué que les étudiant\m es observé\m es portent une une attention forte à la démarche mise en \oe uvre dans la preuve, et valorisent une bonne idée, même mal exploitée. 

Dans la section suivante, nous complétons notre analyse à partir des évaluations réalisées individuellement en devoir maison, à la lumière de l'étude du travail des binômes observés.

\subsection{Quelques résultats sur les évaluations réalisées à la maison}

Comme nous l'avons vu dans la partie précédente, la méthodologie d'entretiens d'explicitation nous a permis d'identifier certaines activités d'étudiant\m es lors d'une évaluation entre pairs. Cependant, nous ne savons pas à quel point ces activités apparaissent chez les étudiant\m es ayant évalué seul\m es. Nous cherchons donc ici des traces de ces activités dans les évaluations en devoir maison.

Il est plus difficile, seulement avec les copies évaluées, d'observer une invalidation : les étudiant\m es barrent rarement les écrits de leur camarade, et n'indiquent pas toujours clairement quand une preuve leur paraît invalide. Pour les validations, on observe des symboles de coches ($\checkmark$), et parfois le\m a correcteur\m rice commente la preuve validée. Ces données donnent peu de traces de l'activité de validation ou d'invalidation des étudiant\m es, notamment les arguments et ressources utilisées. 

\begin{figure}
    \centering
    \includegraphics[width=\textwidth]{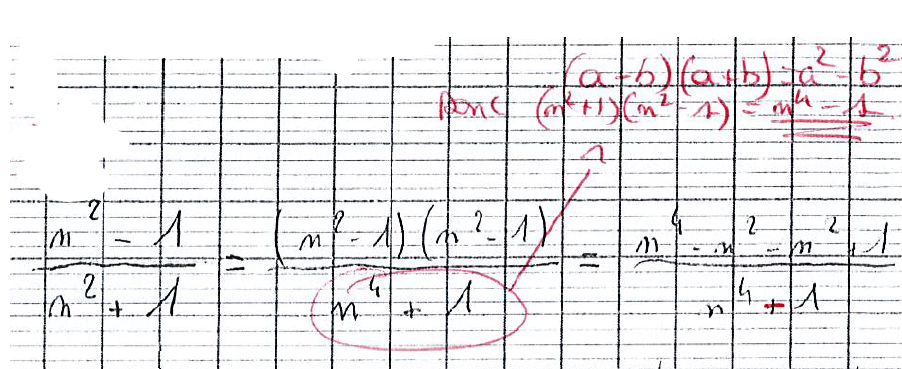}
    \caption{Exemple de correction d'un calcul}
    \label{cop30_calcul}
\end{figure}

\begin{figure}
    \centering
    \includegraphics[width=\textwidth]{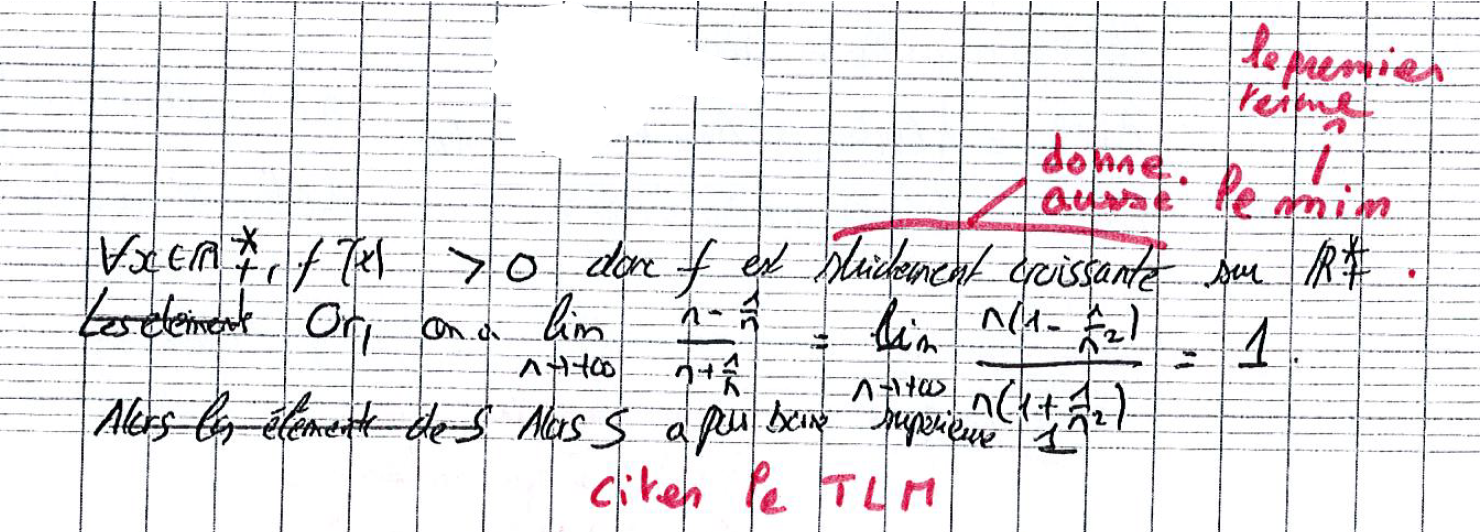}
    \caption{Exemple de complétion d'une preuve}
    \label{cop19_preuve}
\end{figure}

\begin{figure}
    \centering
    \includegraphics[width=\textwidth]{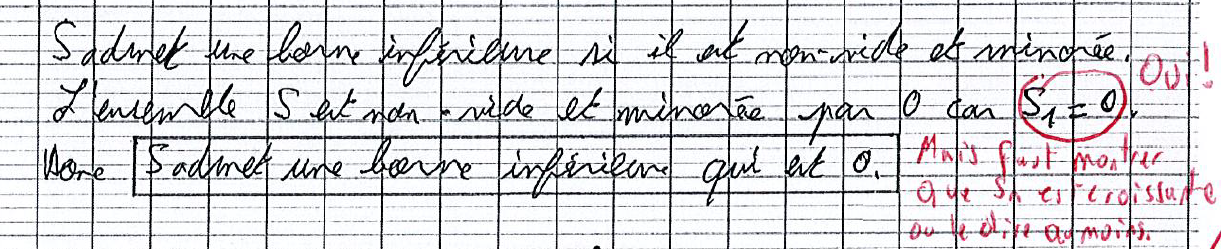}
    \caption{Exemple de proposition d'une stratégie alternative}
    \label{cop29_strat}
\end{figure}

\begin{figure}
    \centering
    \includegraphics[width=0.4\textwidth]{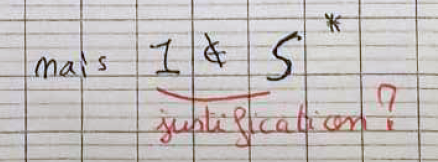}
    \caption{Exemple d'identification de justification manquante}
    \label{cop24_justif}
\end{figure}

Nous avons par contre accès aux corrections proposées par l'évaluateur\m rice. On retrouve des motifs identifiés dans la partie précédente : correction de calculs ou de définitions (figure \ref{cop30_calcul}), complétion d'une preuve (figure \ref{cop19_preuve}), proposition d'une autre stratégie (figure \ref{cop29_strat}), identification d'une justification manquante (figure \ref{cop24_justif}).

\begin{figure}
    \centering
    \includegraphics[width=\textwidth]{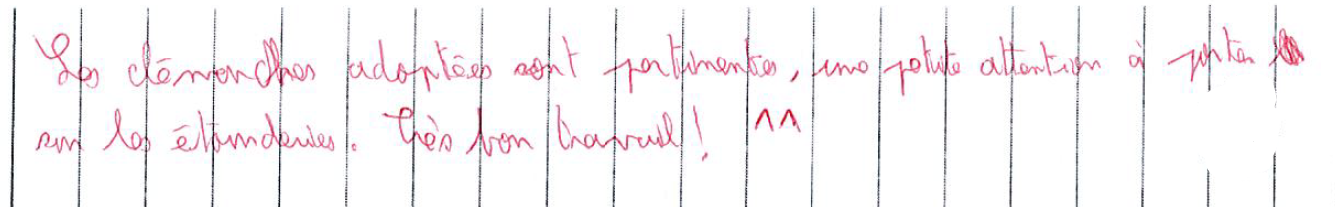}
    \caption{Exemple d'un correcteur qui souligne la pertinence des démarches}
    \label{cop3_demarche}
\end{figure}

Comme nous l'avons rapidement noté dans la section précédentes, les étudiant\m es prêtent attention à la démarche entamée dans la copie qu'ils évaluent. Plusieurs commentaires dans les copies évaluées à la maison vont aussi dans ce sens (figure \ref{cop3_demarche}). Nous avons également remarqué de nombreux commentaires sur la propreté des copies, ainsi que des conseils pour rendre la rédaction plus fluide.

Notons enfin que certains annotations des copies tendent à montrer que les étudiant\m es ont utilisé des références externes, les éléments de correction à disposition, notamment dans le second exercice, moins maîtrisé que l'exercice 1. 



\subsection{Synthèse des analyses de l'activité des étudiant\m es}

Les analyses des entretiens et des copies nous ont permis d'identifier les activités des étudiant\m es en situation d'évaluer des pairs. Les trois axes principaux sont la compréhension des preuves, la validation ou l'invalidation de ces preuves, et leur correction. En plus d'évaluer le fond du contenu mathématique, les étudiant\m es sont attentif\m ves aux méthodes de résolution employées et au résultat final obtenu. Chez certain\m es, ces observations sur la démarche peuvent servir à l'invalidation de la preuve évaluée, mais chez d'autres la démarche est seulement commentée.

Nous avons également observé que les étudiant\m es relèvent les imprécisions de syntaxe et de rédaction, cependant il est difficile de comprendre quel est l'impact sur la validation ou l'invalidation des preuves.

Notons que la lecture, la compréhension et la validation des preuves s'effectuent en général en parallèle. Pour ce faire, les étudiant\m es décomposent souvent la preuve en ses grands axes et éléments importants, pour essayer d'en comprendre la structure générale. Pendant les moments de validation, les étudiant\m es s'interrogent régulièrement sur les justifications nécessaires à la validité d'une preuve ou d'un calcul. Ces réflexions pourront être utiles à une compréhension plus en profondeur de l'activité de preuve et des concepts en jeu dans les exercices.

On remarque que la vigilance des étudiant\m es aux erreurs dans les preuves varie. Il est arrivé plusieurs fois que les étudiant\m es, pensant reconnaître une méthode classique ou un théorème du cours, ne fassent pas attention à la présence de toutes les hypothèses. Ainsi, certaines copies qui utilisent le théorème de la limite monotone sans vérifier la majoration de la suite en jeu ne sont pas invalidées par les pairs. On trouve aussi des copies utilisant le théorème de la limite monotone sans le citer explicitement mais en vérifiant toutes les hypothèses, et qui sont invalidées. On constate ici un manque de recul de certain\m es étudiant\m es sur la validation de l'utilisation d'un théorème. Les étudiant\m es semblent aussi plus critiques envers des méthodes qu'il\m elles ne connaissent pas, ou n'ont pas eux-mêmes ou elles-mêmes mises en place. 


Cette catégorisation de l'activité des étudiant\m es nous permet donc de préciser les éléments qu'il\m elles prennent en compte lors de l'évaluation de pairs. Ils peuvent concerner les preuves en elles-mêmes, les résultats présentés, la démarche entreprise ou la rédaction et la syntaxe. Les étudiant\m es sont aussi attentif\m ves aux éléments non mathématiques, comme la propreté ou la clarté de la copie qu'il\m elles évaluent. Certain\m es s'intéressent également à la qualité et l'efficacité des méthodes empruntées. Tout ceci est sans doute lié au contrat didactique de la classe vis-à-vis de l'évaluation, très spécifique en CPGE.

Pour compléter nos analyses de l'activité des étudiant\m es, nous avons mis en place un questionnaire qui a permis de recueillir la perception qu'ont les étudiant\m es de leur activité. Nous synthétisons les résultats dans la section suivante.

\section[Perception des étudiant\m es sur leurs activités]{Perception des étudiant\m es sur leurs activités : retours du questionnaire}
\label{sectionquestionnaire}

Comment les étudiant\m es ont-il\m elles perçu l'activité d'évaluation entre pairs à laquelle il\m elles ont participé ? Ont-il\m elles pris conscience de certaines choses au cours de ce processus (relativement à la preuve, à l'évaluation, à la rédaction, au contrat didactique\dots) ? Pour étudier cette question, nous avons fait passé un court questionnaire aux étudiant\m es. 

Il a été transmis à la fin de l'expérimentation, et a reçu 38 réponses. Les premières questions se concentrent sur le ressenti des étudiant\m es quant aux apports de cette expérimentation, de manière générale et sur la compréhension des exercices. Nous questionnons ensuite les critères d'évaluation qui ont semblé importants aux étudiant\m es. Les dernières questions interrogent les étudiant\m es sur leurs impressions sur la posture d'évaluateur\m rice.

Le détail des questions et un tableau résumant les statistiques des réponses aux questions fermées sont disponibles en annexe \ref{Qtaire}. Nous résumons ci-dessous les observations qui ressortent des réponses au questionnaire, relativement à nos questions de recherche. Il faut noter que ces résultats émanent de déclarations des étudiant\m es, et attestent donc de la perception qu'il\m elles ont de l'activité proposée, qui peut différer de leur activité réelle. Il est néanmoins intéressant de savoir comment il\m elles ont perçu leur activité et, sur le plan "méta", ce dont il\m elles ont pu prendre conscience lors de l'évaluation entre pairs.

\paragraph{Critères d'évaluation jugés importants}

Nous avons demandé aux étudiant\m es les critères d'évaluations qu'il\m elles estiment être les plus importants, qui sont synthétisés ci-dessous :

\begin{itemize}
    \item la clarté de la rédaction (12 étudiant\m es);
    \item la rigueur dans les justifications et l'utilisation d'un théorème (9 étudiant\m es);
    \item la propreté de la copie (7 étudiant\m es);
    \item la démarche de résolution (6 étudiant\m es);
    \item le résultat final (5 étudiant\m es);
    \item la connaissance précise du cours (3 étudiant\m es);
    \item le fait de répondre aux questions posées (2 étudiant\m es);
    \item le fait de ne pas bluffer (1 étudiant\m e).
\end{itemize}

Notons que les aspects de communication du raisonnement sont les plus représentés dans les réponses (clarté de la rédaction et propreté de la copie), ce qui semble être une prise de conscience pour une partie des étudiant\m es. Des critères de validité mathématique sont aussi évoqués par les étudiant\m es, qui insistent sur l'importance de la rigueur dans les justifications et la démarche de résolution. Nous avons également relevé que 9 étudiant\m es identifient le critère de "l'importance du raisonnement" comme important, sans que l'on puisse tirer plus de détails quant aux caractéristiques du raisonnement auxquelles il\m elles se référent.

\paragraph{Perception sur les apprentissages mathématiques}

Les réponses des étudiant\m es semblent indiquer qu'évaluer des pairs est perçu comme contribuant effectivement à leur apprentissage. Cependant, dans leur perception, être soi-même évalué par des pairs s'avère moins utile à la compréhension des exercices, et est très dépendant de la personne qui évalue.

Des étudiant\m es notent qu'évaluer des pairs nécessite une maîtrise avancée du cours. De plus, il\m elles trouvent utile d'être confronté\m es à plusieurs résolutions d'un même exercice, ce qui permet \enquote{d'échanger des savoirs et des savoir faire}. De plus, certain\m es étudiant\m es notent que l'évaluation d'autres copies les a amené\m es à confronter leurs propres réponses à celles qu'il\m elles ont évaluées, ce qui leur permet de \enquote{prendre du recul sur [leur] propre travail}. Nous remarquons ici que ces réponses font échos aux références aux expériences personnelles (utilisées comme références externes) mobilisées par les étudiant\m es lors des entretiens.

Les étudiant\m es ont également relevé l'intérêt de rencontrer plusieurs erreurs de résolution, en particulier des erreurs qu'il\m elles auraient pu faire aussi.

Plus largement, certain\m es affirment que l'évaluation de copies de pairs leur a permis de mieux saisir les "étapes clés" qui constituent un raisonnement. Cela leur permet donc de comprendre la structure générale d'une preuve qu'ils rencontrent. On peut supposer que les étudiant\m es perçoivent un intérêt qui va au-delà de l'exercice spécifique traité, sur la preuve par exemple.

\paragraph{Perception des aspects méta-cognitifs}

Comme nous l'avons cité ci-dessus, une partie des étudiant\m es disent se rendre mieux compte de l'importance d'une rédaction claire et soignée, donc de l'importance de la communication de la preuve en plus de la preuve en elle-même. Par ailleurs dans le questionnaire, tout comme dans les entretiens, plusieurs étudiant\m es ont remarqué que les "arnaques" peuvent être facilement identifiées par le correcteur\m  rice.

On peut interpréter tout cela comme relevant d'une prise de conscience de la dimension humaine de l'activité de preuve, et de la nécessité de tenir compte du point de vue de lecteur\m rice lors de la production d'une preuve.

Les étudiant\m es relèvent donc plusieurs raisons pour lesquels le processus d'évaluation par les pairs leur a semblé utile. Cela leur permet d'examiner d'autres rédactions pour leur propre stratégie de résolution, et de découvrir d'autres solutions à un même exercice. De plus, cela leur donne l'occasion de réfléchir à l'importance de la structuration d'une preuve dans sa validité. Enfin, il\m elles développent un autre regard sur l'activité d'évaluation et sur la communication d'un contenu mathématique.

\section{Conclusions}

\subsection{Résultats et perspectives}

Nous avons expérimenté un processus d'évaluation entre pairs avec des étudiant\m es du début de l'enseignement supérieur, qui s'est avéré fonctionnel (les étudiant\m es sont bien entré\m es dans la tâche proposée) et qui a engendré une activité riche chez les étudiant\m es.

Nos questions de recherches portaient sur l'activité des étudiant\m es évaluateur\m rices pendant une évaluation de pairs (sur les plans mathématique et cognitif, mais aussi méta-mathématique et méta-cognitif), et plus particulièrement les éléments qu'il\m elles prennent en compte lors des évaluations, les corrections qu'il\m elles apportent, ainsi que les ressources auxquelles il\m elles font appel.

Pour répondre à ces questions, après avoir sélectionné deux exercices portant sur les bornes supérieure et inférieure et sur les manipulations d'ensembles, nous avons réalisé une analyse a priori des différentes stratégies de résolution des exercices, puis une analyse a priori des évaluations possibles de ces exercices. Nous avons ensuite observé des étudiant\m es évaluer des copies, afin d'identifier leur activité. Enfin, sur la base de ces analyses, nous avons examiné les traces écrites des évaluations réalisées en devoir maison.

\paragraph{Activités des étudiant\m es lors de l'évaluation}

Lors de l'évaluation d'une copie, les étudiant\m es se trouvent dans des phases de lecture et de compréhension de preuves, des phases de validation et d'invalidation et des phases de correction.

Lors de la lecture et la compréhension de la copie, les évaluateur\m rices cherchent à expliquer et expliciter le raisonnement entrepris dans la copie, notamment en identifiant les différentes étapes qui constituent la preuve. Pour cela, il\m elles peuvent faire référence aux contenus de leur cours afin de préciser les concepts et énoncés mathématiques en jeu. Parfois, il\m elles sont amené\m es à interpréter ce que la personne évalué\m e "a voulu faire" dans la copie. 

Lors de la phase de validation et d'invalidation, les étudiant\m es peuvent s'intéresser à la preuve dans son ensemble, à une étape du raisonnement, aux propriétés des objets manipulés, au résultat ou encore à la démarche de la preuve. Pour justifier une invalidation, il\m elles peuvent être amené\m es à trouver des contre-exemples à des énoncés mathématiques utilisés, implicitement ou explicitement, dans la copie. Nous avons observé ceci particulièrement lorsque l'étudiant\m e évalué\m e semble utiliser un théorème présent dans le cours (comme le théorème de la limite monotone), sans en vérifier toutes les hypothèses. Enfin, les évaluateur\m rices discutent de la précision nécessaire dans les justifications permettant de valider une preuve, et identifient les résultats "pas assez justifiés".

Nous notons de plus qu'au cours de l'évaluation, une attention particulière est portée par les étudiant\m es à la démarche de résolution de l'exercice empruntée dans la copie. Si ce n'est pas forcément un critère de validation de la preuve, de nombreux étudiant\m es commentent lorsqu'il\m elles estiment que la personne corrigée a eu "de bonnes idées" dans la preuve. Les étudiant\m es prennent également en compte le résultat final de la preuve, dont la validité peut influencer leur opinion de correcteur\m rices sur la validité du reste de la preuve.

Enfin, les étudiant\m es proposent des corrections des éléments qu'il\m elles ont invalidés. Il peut s'agir d'erreurs de définition, qui sont systématiquement corrigées, ou d'erreurs de calcul, pour lesquelles les étudiant\m es se contentent parfois de noter que le résultat est erroné. Pour la correction des éléments relatifs à la preuve, les évaluateur\m rices peuvent compléter le raisonnement entamé dans la copie ; nous avons alors observé des discussions portant sur le caractère nécessaire de certaines arguments. Certain\m es étudiant\m es peuvent aussi changer la structure d'une preuve en modifiant l'ordre des différentes étapes de raisonnement, pour la rendre valide. Enfin, il arrive que des étudiant\m es proposent une stratégie de résolution différente de celle entamée dans la copie.

Ainsi, on voit que le processus d'évaluation entre pairs peut amener les étudiant\m es à une activité très riche autour de l'analyse de solutions et de preuves, impliquant certaines réflexions que l'on peut considérer d'ordre méta-mathématique.

\paragraph{Utilisation des arguments et références internes et externes}

Les arguments et références mobilisées par les étudiant\m es lors de l'évaluation n'ont pas seulement une visée de validation : ils peuvent aussi servir à éclairer des preuves qui n'ont pas été initialement comprises par les évaluateur\m rices, ou sont utilisés pour commenter les raisonnements proposés dans la copie. Le principal argument interne que nous avons relevé est la production de contre-exemple pour invalider un énoncé affirmé dans la copie. D'autres arguments internes, références au cours et vérification d'éléments techniques (formules, calculs\dots) ont également été observés. Nous avons noté peu de références explicites au cours ou à des sources expertes lors des entretiens. De nombreux étudiant\m es ont fait référence à leur propre résolution des exercices, soit pour commenter les preuves évaluées, soit pour valider ou invalider un calcul ou le résultat d'une étape du raisonnement. Enfin, la feuille de correction que nous avons fournie a été exploitée à plusieurs reprises lors des entretiens, et semble avoir été également utilisée pour la correction de copies en devoir maison.

En conclusion, il semble que les étudiant\m es mobilisent à la fois des arguments internes et externes, avec une variabilité selon les individus. La dialectique entre ces deux dimensions dans les processus d'évaluation de pairs, et le développement des arguments internes chez les étudiant\m es seraient certainement intéressants à étudier sur un temps plus long.

\paragraph{Apports de l'évaluation entre pairs à l'apprentissage}

Les analyses des entretiens, des copies et du questionnaire réalisées nous ont permis de dégager des aspects de l'évaluation entre pairs pouvant participer à l'apprentissage des étudiant\m es. Ceci va dans le sens des résultats des recherches existantes mentionnées en section \ref{generalites}.

Notre étude montre que ce processus peut permettre de développer un travail sur les contenus mathématiques en jeu. En effet, l'évaluation de copies de pairs mets les étudiant\m es face à d'autres stratégies de résolution que celles qu'il\m elles auraient entrepris. Cela peut permettre de se rendre compte de la diversité de moyens possibles pour résoudre un problème mathématique. À l'inverse, les étudiant\m es peuvent évaluer des preuves similaires à celles qu'il\m elles ont mises en place. Il\m elles ont alors l'occasion de considérer à nouveau à la validité de leur méthode, et d'observer d'autres erreurs qui peuvent survenir.

De plus, l'évaluation et la correction d'éléments techniques peut permettre aux étudiant\m es de travailler le cours. Par exemple, lorsqu'il\m elles cherchent des contre-exemples à certaines étapes de déduction, il\m elles sont amené\m es à réfléchir à l'importance des différentes hypothèses des théorèmes qui sont utilisés. On peut s'attendre à ce que cela permette une compréhension plus profonde du cours et des théorèmes mis en jeu, qui pourrait mener, plus généralement, à saisir l'importance des hypothèses dans les théorèmes.

L'évaluation par les pairs est également utile pour prendre du recul sur l'activité de preuve. Lors des différentes phases de l'évaluation décrites ci-dessus, les étudiant\m es peuvent s'interroger sur la structure des preuves qu'il\m elles évaluent, mettant en lumière l'importance des étapes d'une preuve et de leur organisation.

Enfin, l'évaluation entre pairs permet d'expliciter certains enjeux de l'évaluation, en rendant explicite le contrat didactique en jeu. Nous avons noté que les étudiant\m es semble mieux appréhender les évaluations d'un\m e correcteur\m rice. En particulier, ce processus semble mettre au jour des enjeux de communication des contenus mathématiques, autant dans la présentation d'un écrit que dans la rédaction des raisonnements.

En conclusion, nous observons que les étudiant\m es dans un processus d'évaluation entre pairs développent une activité mathématique qui diffère de celle que l'on peut imaginer ordinairement dans la résolution des mêmes exercices. Nous observons une activité d'analyse de preuve, d'interprétation des raisonnements, d'évaluation de la validité d'étapes ou de l'entièreté d'une preuve, ainsi que le développement des réflexions sur les plans méta-mathématiques et méta-cognitifs.
Ceci nous semble porteur d'un potentiel pour l'apprentissage mathématiques, en particulier au début de l'enseignement supérieur.

Nous identifions certaines limites à ce travail : sur le plan de la recherche, notre travail exploratoire montre un potentiel pour l'apprentissage des mathématiques, mais il reste à comprendre comment peut se réaliser ce potentiel et le confirmer dans d'autres contextes (scolaires et mathématiques), et à étudier sur le long terme les effets d'un tel dispositif sur les apprentissages. Sur le plan du dispositif lui-même, l'expérience peu fructueuse du second exercice montre que les conditions de fonctionnement d'une évaluation entre pairs sont importantes et à bien contrôler.

Ce travail ouvre des perspectives de recherches didactiques sur les mathématiques dans le supérieur. Une étude plus approfondie serait à développer sur l'activité des étudiant\m es en situation d'évaluation de pairs, avec un corpus contenant d'autres tâches et dans des contextes divers, ainsi que sur des mises en \oe uvres plus longues et répétées. Un approfondissement des activités possibles des étudiant\m es en situation d'évaluation de pairs, en particulier en croisant avec les recherches sur l'apprentissage de la preuve serait certainement fructueux. Des recherches sur la viabilité en classes de mathématiques de dispositifs d'évaluation entre pairs seraient aussi intéressantes à déployer, avec de comprendre les conditions et contraintes de fonctionnement de ces dispositifs.

\subsection{Quelques conseils et pistes de mise en place}

Sur la base de nos lectures et de notre expérience pour cet article, nous formulons quelques idées, pistes, conseils pour les enseignant\m es qui souhaiteraient mettre en place une évaluation entre pairs en mathématiques, en particulier au début de l'enseignement supérieur. Ces recommandations et points de vigilance ne sont pas tous à prendre comme des résultats étayés, mais il nous a semblé intéressant, dans cette section distincte des résultats de notre étude, de proposer quelques points qui peuvent aider à développer l'évaluation entre pairs en mathématiques.

Au préalable de l'expérimentation, lors du choix du ou des exercices à faire évaluer, il nous semble important :
\begin{itemize}
    \item que les étudiant\m es soient à l'aise avec le type d'exercice, voire que celui-ci porte sur des connaissances suffisamment anciennes et stabilisées ;
    \item qu'une partie de l'exercice porte sur la production d'une preuve (pas seulement d'un calcul) ;
    \item que l'exercice puisse être résolu avec plusieurs stratégies, et que plusieurs points de vue soient possibles sur les objets en jeu (en ce sens, un énoncé "ouvert" qui ne donne pas la réponse, c'est-à-dire qui ne soit pas de la forme "montrer que", nous semble assez adapté).
\end{itemize}
 
Pour la mise en place du processus d'évaluation entre pairs lui-même, nous recommandons de :

\begin{itemize}
    \item clarifier les attentes et les objectifs aux étudiant\m es (par exemple faire un lien avec l'activité d'auto-évaluation), et préciser si ce processus permettra de donner une note en tant qu'évaluateur\m rice et/ou en tant que qu'évalué\m e ; 
    \item développer et clarifier les critères d'évaluation en impliquant les étudiant\m es, et en particulier discuter des critères de validation d'une preuve possibles (références internes et externes) -- On peut à cette occasion donner des exemples d'évaluation ;
    \item au cours des premières sessions, guider, encadrer et modérer les corrections et feedback donnés par les évaluateur\m rices.
    \item faire un retour aux évaluateur\m rices après le processus, surtout les premières fois.
\end{itemize}

Nous suggérons également que tous\m tes les étudiant\m es cherchent au préalable l'exercice qu'il\m elles vont évaluer, et soient tous\m tes évalué\m es par leurs pairs sur cet exercice.

\begin{figure}[h!]
    \centering
    \includegraphics[width = 0.45\textwidth]{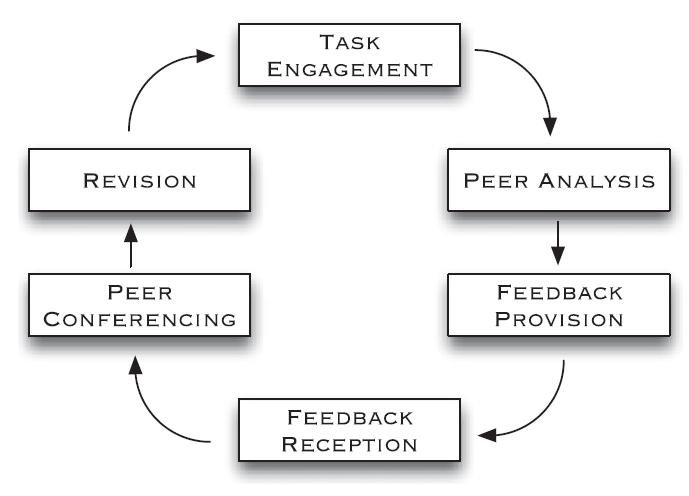}
    \caption{Le cycle d'évaluation de \citet{reinholz_assessment_2015}}
    \label{asscycle}
\end{figure}

Il peut être intéressant de mettre en place un "cycle d'évaluation" (assessment cycle) tel que le propose \citet{reinholz_assessment_2015} (figure \ref{asscycle}).

Nous pointons enfin le fait que ce processus d'évaluation nécessite que les étudiant\m es aient déjà une idée du contrat didactique concernant l'évaluation dans la formation concernée. Ainsi, il faut certainement un temps d'appropriation d'un nouveau contrat didactique d'évaluation avant d'initier des évaluations entre pairs. Les classes préparatoires aux grandes écoles nous ont semblé être appropriées pour l'expérimentation, et nous nous questionnons sur les conditions de mise en place de ce processus à l'université, certainement plus complexes.


Il est souvent recommandé de développer l'évaluation entre pairs sur plusieurs sessions au fil d'un semestre ou d'une année. Cependant, selon notre expérience, même une seule session d'évaluation entre pairs semble initier une prise de conscience de certains enjeux aux étudiant\m es (enjeux de preuve, enjeux d'évaluation\dots) et leur faire se poser des questions formatrices. Nous encourageons donc nos collègues à expérimenter, même ponctuellement, un dispositif d'évaluation entre pairs. 

\begin{appendices}

\section{Grille d'analyse}
\label{annexegrillenalayse}

\begin{tabular}{|c|p{6cm}|p{9cm}|}
\hline

\multicolumn{3}{|c|}{\textit{Évaluations portant sur le fond du raisonnement mathématique}}\\
\hline

\textbf{1}& \multicolumn{2}{c|}{\textbf{Activités portant sur l'évaluation d'une preuve proposée dans la copie}}\\
\hline
1Z & Lire la preuve & \\
\hline

1A & Essayer de comprendre la preuve & Identifier des étapes du raisonnement mathématique \\
 & & Reconnaître les concepts mathématiques utilisés \\
 & & Expliquer une étape du raisonnement \\
 & & Anticiper un raisonnement \\
\hline
1B & Évaluer la validité de la preuve & Invalider / Valider la démarche \\
 & & Invalider / Valider  une étape du raisonnement / d'un calcul \\
 & & Invalider / Valider la preuve \\
\hdashline
 & & Noter l'absence de la réponse à une question \\
 & & Pointer un manque de justification ou une bonne justification \\
 & & Pointer une information non nécessaire au raisonnement \\

\hline
1C & Évaluer la vraisemblance de la preuve & Évaluer la vraisemblance de la démarche \\
 & & Évaluer la vraisemblance d'une étape / d'un calcul \\
 & & Évaluer la vraisemblance de la preuve \\

\hline

1D & Corriger la preuve & Corriger une étape du raisonnement \\
 & & Corriger la structure d'une preuve \\
 & & Compléter un raisonnement \\
 & & Proposer une autre méthode de résolution \\
 & & Corriger un détail technique \\
\hdashline
 & & Apporter des précisions sur une étape du raisonnement \\
 & & Expliquer les raisons d'une erreur \\
\hline
\hline

\textbf{2}& \multicolumn{2}{c|}{\textbf{Activités portant sur l'évaluation d'une preuve envisagée par l'élève corrigé}}\\
\hline

2A & Essayer de comprendre la preuve & \\
\hline
2B & Évaluer la validité de la preuve & \\
\hline
2C & Évaluer la vraisemblance de la preuve &  \\
\hline
2D & Corriger la preuve & \\
\hline
\hline

\textbf{3}& \multicolumn{2}{c|}{\textbf{Activités pourtant sur l'évaluation du résultat d'une preuve}}\\
\hline

3B & Évaluer la validité du résultat & \\
\hline
3D & Corriger le résultat & \\
\hline

\end{tabular}
\\
\begin{tabular}{|c|p{6cm}|p{9cm}|}
\hline

\multicolumn{3}{|c|}{\textit{Évaluations portant sur la forme du raisonnement mathématique}}\\
\hline

\textbf{4}& \multicolumn{2}{c|}{\textbf{Activités pourtant sur l'évaluation de la rédaction et de la syntaxe d'une preuve}}\\
\hline

4B1 & Évaluer la validité de la rédaction et de la syntaxe &  Valider / Invalider une définition \\
 &  &  Valider / Invalider la syntaxe d'un énoncé \\
\hline
4B2 & Évaluer la clarté de la rédaction & \\
\hline
4D & Corriger la rédaction ou la syntaxe & \\

\hline
\end{tabular}

\begin{tabular}{|c|p{6cm}|p{9cm}|}
\hline
\multicolumn{3}{|c|}{\textit{Évaluations ne portant pas sur le raisonnement mathématique}}\\
\hline

\textbf{5}& \multicolumn{2}{c|}{\textbf{Activités ne portant pas sur l'évaluation du raisonnement mathématique}}\\
\hline
5A & Juger de la propreté de la copie & \\
\hline
5B & Essayer de comprendre le comportement de l'élève & \\
\hline
5C & Juger de la qualité ou de la légitimité d'une méthode &  \\
\hline
5D & Juger de la compréhension d'un élève & \\
\hline

\end{tabular}

\section{Questionnaire}
\label{Qtaire}

Le questionnaire a été partagé aux étudiant\m es via la plateforme Framaforms. Les questions sont les suivantes :

\begin{enumerate}
    \item Nom
    \item Prénom
    \item Étiez-vous volontaire pour faire les corrections sous ma supervision ?
    \item Le fait de corriger des copies d'autres élèves vous a-t-il semblé utile ? Pourquoi ?
    \item Le fait d'être corrigé par d'autres élèves vous a-t-il été utile ? Pourquoi ?
    \item Est ce que le fait de corriger des copies d'autres élèves vous a permis de mieux comprendre les exercices ? Pourquoi ?
    \item Est ce que le fait d'être corrigé par d'autres élèves vous a permis de mieux comprendre les exercices ? Pourquoi ?
    \item Quels sont les critères d'évaluation les plus importants que vous retenez ?
    \item Qu'est ce qui a été le plus facile pour vous dans le fait de corriger des copies ?
    \item Qu'est ce qui a été le plus difficile pour vous dans le fait de corriger des copies ?
    \item Vous êtes-vous senti à l'aise en corrigeant les copies ? Pourquoi ?
    \item Est-ce une expérience que vous aimeriez refaire ? Pourquoi ?
\end{enumerate}

Le tableau ci-dessous synthétise les réponses aux questions fermées. 38 étudiant\m es ont répondu à ce questionnaire.

\begin{tabular}{|p{8cm}|c|c|}
\hline 
     Question & Oui et plutôt oui & Non et plutôt non \\
     & Utile et plutôt utile & Inutile et plutôt inutile \\
\hline
    Le fait de corriger des copies d'autres élèves vous a-t-il semblé utile ? & 36 & 2\\
\hline
    Le fait d'être corrigé par d'autres élèves vous a-t-il été utile ? & 27& 11\\
\hline
    Est ce que le fait de corriger des copies d'autres élèves vous a permis de mieux comprendre les exercices ? & 27&11\\
\hline
    Est ce que le fait d'être corrigé par d'autres élèves vous a permis de mieux comprendre les exercices ? & 17&21\\
\hline
    Vous êtes-vous senti à l'aise en corrigeant les copies ? & 26&12\\
\hline
    Est-ce une expérience que vous aimeriez refaire ? & 34&4\\
\hline
\end{tabular}
\end{appendices}

\bibliographystyle{apacite}
\bibliography{epidemes_references}

\authoraddresses{
Juliette Veuillez-Mainard \\
Institut Montpellierain Alexander Grothendieck, Université de Montpellier, CNRS\\
\email juliette.veuillez-mainard@umontpellier.fr

Simon Modeste \\
Institut Montpellierain Alexander Grothendieck, Université de Montpellier, CNRS\\
\email simon.modeste@umontpellier.fr
}

\end{document}